%% file: fw-gentle.tex
\PassOptionsToPackage{dvipsnames}{xcolor}
\documentclass{article}
\newif\ifedit
 \editfalse

\input{prologue.tex}

\usepackage{todonotes}

\usepackage{tabstackengine}
\setstackgap{L}{1.2\normalbaselineskip}
\setstacktabbedgap{.4em}
\fixTABwidth{T}
\begin{document}

\title{The Frank-Wolfe algorithm: a short introduction}

\author{\name Sebastian Pokutta \email \href{mailto:pokutta@zib.de}{pokutta@zib.de} \\
       \addr 
%       Institute of Mathematics \& AI in Society, Science, and Technology\\
       Zuse Institute Berlin \& TU Berlin \\
       Berlin,  Germany}

\maketitle
  
\begin{abstract}
In this paper we provide an introduction to the Frank-Wolfe algorithm, a method for smooth convex optimization in the presence of (relatively) complicated constraints. We will present the algorithm, introduce key concepts, and establish important baseline results, such as e.g., primal and dual convergence. We will also discuss some of its properties, present a new adaptive step-size strategy as well as applications. 
\end{abstract}

\medskip \hfill (Article for the \emph{Jahresbericht der Deutschen Mathematiker Vereinigung}) \medskip

\section{Introduction\pagehint{2}}\label{sec:introduction} 

Throughout this paper we will be concerned with \emph{constrained optimization problems} of the form
\begin{equation}
	\label{eq:mainProb}
	\tag{Opt}
	\min_{x \in P} f(x),
\end{equation}
where $P \subseteq \RR^n$ is some convex feasible region capturing the constraints, e.g., a polyhedron arising from a system of linear inequalities or a spectrahedron, and $f$ is the objective function satisfying some regularity property, e.g., smoothness and convexity. We also need to specify what access methods we have, both, to the function and the feasible region. A common setup is black box first-order access for $f$, allowing (only) the computation of gradients $\nabla f(x)$ for a given point $x$ as well as the function value $f(x)$. For the access to the feasible region $P$, which we will assume to be compact in the following, there are several common models; we simplify the exposition here for the sake of brevity:
\begin{enumerate}% [left=0pt]
	\item \emph{Projection.} Access to the projection operator $\Pi_P$ of $P$ that, for a given point $x \in \R^n$ returns $\Pi_P(x) \doteq \argmin_{y \in P} \norm{x-y}$, for some norm $\norm{.}$ (or more generally Bregman divergences).
	\item \emph{Barrier function.} Access to a barrier function of the feasible region $P$ that increases in value to infinity when approaching the boundary of $P$. A typical example is, the barrier function $- \sum_i \log (b_i - A_i x)$ for a linear inequality system $P \doteq \{x \mid Ax \leq b\}$.
	\item \emph{Linear Minimization.} Access to a Linear Minimization Oracle (LMO) that, given a linear objective $c \in \R^n$, returns $y \in \argmin_{x \in P} \innp{c}{x}$.
\end{enumerate}
Specialized approaches for specific cases, e.g., the simplex method \citep{dantzig1981reminiscences,dantzig1983reminiscences} in the case of linear objectives which uses an explicit description of the feasible region also exist but here we concentrate on the aforementioned black box model. There are also proximal methods, which can be considered generalizations of projection-based methods and which we will not explicitly consider for the sake of brevity; see e.g., \citet{nemirovsky1983problem,nesterov04,nesterov18,Nocedal} for a discussion.	

\begin{figure}[b]
	\begin{minipage}[b]{0.3\textwidth}
		\centering
		\includegraphics[scale=0.65]{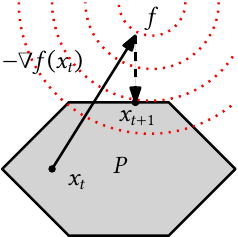}
		
		\caption{\label{fig:pgd}Projection-based methods: may require 
			(potentially expensive) projection back into $P$ to ensure feasibility.}
	\end{minipage} 
	\hfill	
	\begin{minipage}[b]{0.3\textwidth}
		\centering
		\includegraphics[scale=0.65]{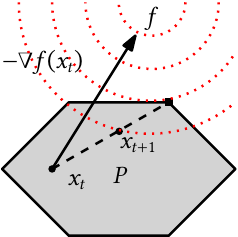}
		
		\caption{\label{fig:fw}Projection-free methods: % typically follow alternative directions and
			ensure feasibility by forming convex combinations only. \\ \ }
	\end{minipage}
	\hfill
	\begin{minipage}[b]{0.35\textwidth}
		\centering
		\begingroup
		%\renewcommand{\ALG@name}{Frank-Wolfe Algorithm}
		%\renewcommand{\thealgorithm}{}
		%\floatname{algorithm}{}
		\begin{algorithm}[H]
			\caption{FW algorithm}
			\begin{algorithmic}[1]
				{\footnotesize
					\STATE$x_0\in P$
					\FOR{$t=0$ \textbf{to} $T-1$}
					\STATE{$v_t\leftarrow\argmin\limits_{v\in P}\langle{\nabla f(x_t)},v\rangle$}
					\STATE{$x_{t+1}\leftarrow (1-\gamma_t) x_t + \gamma_t v_t$}
					\ENDFOR}
			\end{algorithmic}
		\end{algorithm}
		\endgroup
%		\vspace{-1\baselineskip}
		
		\caption{\label{fig:alg} New iterates are formed by convex combination with an extreme point approximating the gradient, ensuring feasibility.}
	\end{minipage}
	
\end{figure}

Traditionally, problems of the form \eqref{eq:mainProb} are solved by variants of \emph{projection-based methods}. In particular first-order methods, such as variants of projected gradient descent are often chosen in large-scale contexts as they are comparatively cheap. For some feasible region $P$ with projector $\Pi_P$ (e.g., $\Pi_P(x) \doteq \argmin_{y \in P} \norm{x-y}$) and smooth objective function $f$, \emph{projected gradient descent (PGD)} updates typically take the form:
\begin{align}
	\tag{PGD}
	x_{t+\nicefrac{1}{2}} & \leftarrow x_t - \gamma_t \nabla f(x_t) \\
	\nonumber 	x_{t+1} & \leftarrow \Pi_P(x_{t+\nicefrac{1}{2}}), 
\end{align}
where $\gamma_t$ is some step-size, e.g., $\gamma_t = 1/L$ if $f$ is
$L$-smooth (see Definition~\ref{def:smooth}) and convex. In essence, a descent step is taken without considering the constraints, and then it is projected back into the feasible region (see Figure~\ref{fig:pgd}). Projection-based first-order methods have been extensively studied, with comprehensive overviews available in, e.g., \citet{nesterov18,Nocedal}. Optimal methods and rates are known for most scenarios. Efficient execution of the projection operation is possible for simple constraints, such as box constraints or highly structured feasible regions, e.g., as discussed in \citet{moondra2021reusing,gupta2016solving} for submodular base polytopes. However, when the feasible region grows in complexity, the projection operation can become the limiting factor. It often demands the solution of an auxiliary optimization problem---known as the \emph{projection problem}---over the same feasible region for \emph{every} descent step. This complexity renders the use of projection-based methods for many significant constrained problems quite challenging; in some cases relaxed projections which essentially compute separating hyperplanes can be used though.

\emph{Interior point methods (IPM)} offer an alternative approach, see e.g., \citet{boyd2004convex,potra2000interior}. To illustrate this approach, consider the goal of minimizing a linear function $c$ over a polytope defined as \(P \doteq \{x \mid Ax \leq b\}\). The typical updates in a path-following IPM resemble:
\begin{align}
	\tag{IPM}
	x_{\mu} & \leftarrow \argmin_{x} \innp{c}{x} + \mu \sum \log (b_i - A_i x),
\end{align}
where the value of $\mu \rightarrow 0$ according to some strategy for $\mu$. Often, these steps are only approximately solved. IPMs, while potent with appealing theoretical guarantees, usually necessitate a barrier function that encapsulates the feasible region's description. In numerous critical scenarios, a concise feasible region description is either unknown or proven to be non-existent. For instance, the matching polytope does not admit small linear programs, neither exact ones \citep{rothvoss14} nor approximate ones \citep{BP2014matching,BP2014matchingJour,sinha2018lower}. Additionally, achieving sufficient accuracy in the IPM step updates often requires second-order information, which can sometimes restrict its applicability. 

Upon closely examining the two methods mentioned earlier, it is clear that both essentially transform the constrained problem~\eqref{eq:mainProb} into an \emph{unconstrained} one. They then either correct updates that violate constraints (as in PGD) or penalize nearing constraint violations (as in IPM). Yet, another category of techniques exists, termed \emph{projection-free methods}, which focus directly on constrained optimization. Unlike their counterparts, these methods sidestep the need for costly projections or penalty strategies and maintain feasibility throughout the process. The most notable variants in this category are the \emph{Frank-Wolfe (FW) methods}---going back to \citet{fw56}---which will be the focus of this article and which are also known as \emph{conditional gradient~(CG) methods}  \citep{polyak66cg}.

Historically, methods like the Frank-Wolfe algorithm garnered limited attention because of certain drawbacks, notably sub-optimal convergence rates. However, there was a notable resurgence in interest around 2013. This revival is largely attributed to shifting requirements and their other, now suddenly relevant properties. Notably, these methods are well suited to handle complicated constraints and possess a low iteration complexity. This makes them very effective in the context of
large-scale machine learning problems (see, e.g.,
\citep{lacoste2013block, jaggi13fw, negiar20, dahik2021robust, jing2023robust}),
image processing (see, e.g., \citep{joulin2014efficient,tang2014co}), quantum
physics (see, e.g., \citep{gilbert1966iterative, DIBKGP2023}),
submodular function maximization (see, e.g.,
\citep{feldman2011unified, vondrak2008optimal, badanidiyuru2014fast, MirzasoleimanBK16, hassani2017gradient, mokhtari2017conditional, submodular2017, submodular2017jour, mokhtari2018decentralized, bach2019submodular}), online learning (see,
e.g.,
\citep{hazan12, zhang2017projection, chen2018projection, garber2020projfree-learn, KRAP2021, zhang2023dynamic}) and
many more  (see, e.g.,
\citep{bolte07lojo, clarkson08, pierucci14, zaid15, wang2016parallel, cheung18, ravi19, hazan2020faster, SelfConcordantFW_2020, CaP2020, MBP2021, CPSZ2021, garber2021NEP-FW, FWjourney2021, WHP2022, chen2023reducing, de2023short, DIBKGP2023, DVP2023, lj16nonconvex}). Moreover, there has been a proliferation of modifications to these methods, addressing many of their historical limitations (see, e.g.,
\citep{freund2017extended, lacoste15, garber2015faster, garber2016linearly, LPZZ2017, BPZ2017, BPZ2017jour, BPTW2018, CP2020boost, TTP2021}) and there is an intricate connection between Frank-Wolfe methods and subgradients methods \citep{bach15dualfw}; see \citet{CGFWSurvey2022} for a
comprehensive exposition. 

Rather than relying on potentially expensive projection operations (see Figure~\ref{fig:fw}), Frank-Wolfe methods use a so-called Linear Minimization Oracle (LMO). This subroutine only involves optimizing a linear function over the feasible region, often proving more cost-effective than traditional projections; see \citet{CP2021} for an in-depth comparison. The nuclear norm ball, along with matrix completion, is a prime example highlighting the difference in complexity. The core updates in Frank-Wolfe methods often rely on the following fundamental update:
\begin{align}
	\tag{FW}
	v_t & \leftarrow \argmin_{v \in P} \innp{\nabla f(x_t)}{v}\\
	\nonumber	x_{t+1} & \leftarrow (1- \gamma_t) x_t + \gamma_t v_t,
\end{align}
where any solution to the $\argmin$ is suitable and $\gamma_t$ follows some step-size strategy, e.g., $\gamma_t = \frac{2}{2+t}$. Essentially, the LMO identifies an alternate direction for descent. Subsequently, convex combinations of points are constructed within the feasible region to maintain feasibility. Viewed through the lens of complexity theory, the Frank-Wolfe methods reduce the optimization of a convex function \(f\) over \(P\) into the repeated optimization of evolving linear functions over \(P\).  A schematic of the most basic variant of the Frank-Wolfe algorithm is shown in Figure~\ref{fig:alg}.

For a more comprehensive exposition complementing this introductory article the interested reader is referred to \citet{CGFWSurvey2022}; the notation has been deliberately chosen to be matching whenever possible. 

\subsection*{Outline} We start with some basic notions and notations in Section~\ref{sec:prelims} and then present the original Frank-Wolfe algorithm along with some motivation in Section~\ref{sec:fw-algo}. We then proceed in Section~\ref{sec:properties} with establishing basic properties, such as e.g., convergence and also provide matching lower bounds. While this is primarily an overview article, we do provide a new adaptive step-size strategy in Section~\ref{sec:adaptive}, which is also available in the \emph{FrankWolfe.jl} julia package. In Section~\ref{sec:applications} we then consider applications of the Frank-Wolfe algorithm and also discuss computational aspects in Section~\ref{sec:computations}.

\section{Preliminaries\pagehint{1}}
\label{sec:prelims}

In the following $\norm{\cdot}$ will denote the $2$-norm if not stated otherwise. Note however that in general other norms are possible and have been used in the context of Frank-Wolfe algorithms. Moreover, for simplicity we assume that $f$ is differentiable, which is a standard assumption in the context of Frank-Wolfe algorithms although non-smooth variants are known (see, e.g., \citet{CGFWSurvey2022} for details). 

For our analysis we will heavily rely on the following key concepts:

\begin{definition}[Convexity and Strong Convexity]
	\label{def:smooth}
	Let $f\colon P \rightarrow \RR$ be a differentiable function. Then $f$ is \emph{convex} if
	\begin{equation}
		\label{convex}
		f(y) - f(x) \geq \innp{\nabla f(x)}{y - x}
		\quad \text{for all } x, y \in P.
	\end{equation}
Moreover, $f$ is \emph{\(\mu\)-strongly convex} if
\begin{equation}
	f(y) - f(x) \geq \innp{\nabla f(x)}{y-x}
	+ \frac{\mu}{2} \norm{y-x}^{2}
	\quad \text{for all } x, y \in P.
	\label{strconvex}
\end{equation}
\end{definition}

\begin{definition}[Smoothness]
	Let $f \colon P \to \RR$ be a differentiable function. Then $f$ is \emph{\(L\)-smooth} if
\begin{equation}
	\label{eq:baseSmooth}
	f(y) - f(x) \leq \innp{\nabla f(x)}{y-x}
	+ \frac{L}{2} \norm{y-x}^{2}
	\quad \text{for all } x, y \in P.
\end{equation}
\end{definition}

\begin{figure}[h]
	% \label{fig:algs}   % todo 
	\begin{minipage}[b]{0.49\textwidth}
		\centering
		\includegraphics[scale=0.55]{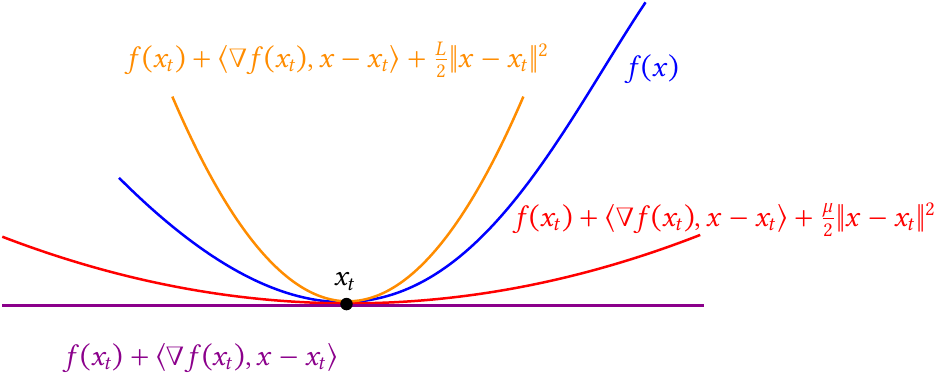}
		
		\caption{\label{fig:conv-smooth} (Strong) convexity and smoothness provide linear and quadratic approximations to $f$. \textcolor{orange}{Orange:} quadratic upper bound via smoothness, \textcolor{red}{Red:} quadratic lower bound via strong convexity, \textcolor{DarkMagenta}{Magenta:} linear lower bound via convexity, \textcolor{blue}{Blue:} function $f$.} 
	\end{minipage} 
	\hfill	
	\begin{minipage}[b]{0.49\textwidth}
		\centering
		\includegraphics[scale=0.55]{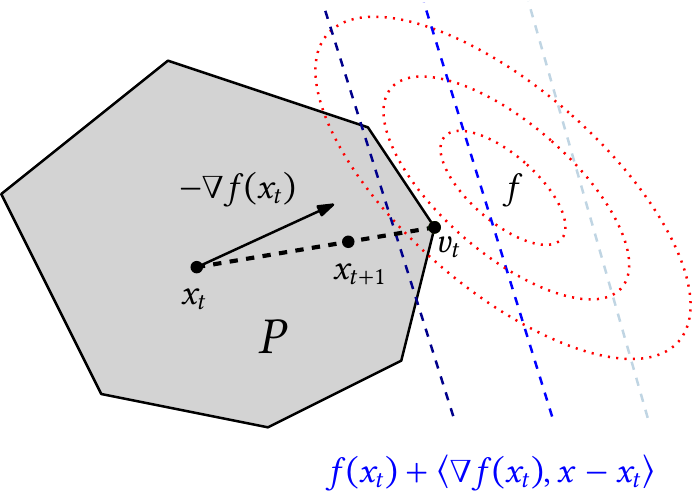}
		
		\caption{
			The Frank–Wolfe step: To minimize a convex function \(f\) over a polytope \(P\), a linear approximation of \(f\) is constructed at \(x_t\) as \(f(x_t) + \innp{\nabla f(x_t)}{x - x_t}\). The Frank–Wolfe vertex \(v_{t}\) minimizes this approximation. The step transitions from \(x_{t}\) to \(x_{t+1}\) by moving towards \(v_{t}\) as determined by a step-size rule. Contour lines of \(f\) in red and linear approximation blue.
			}
		\label{fig:schematic_CG_step}
	\end{minipage}
	
\end{figure}

The smoothness and (strong) convexity inequalities from above allow us to obtain upper and lower bounds on the function $f$. Convexity and strong convexity provide respectively linear and quadratic lower bounds on the function $f$ at a given point $x$ while smoothness provides a quadratic upper bound as shown in Figure~\ref{fig:conv-smooth}.

For completeness we note that, both, $L$-smoothness and $\mu$-strong convexity can also be expressed without relying on function values of $f$ only using gradients $\nabla f$. This is in particular useful in the context of an adaptive step-size strategy that we will discuss in Section~\ref{sec:adaptive} as it significantly improves numerical stability of the estimates. 

\begin{remark}[Smoothness and Strong Convexity via Gradients]
\label{rem:smoothSC}
Let $f\colon P \rightarrow \RR$ be a differentiable function. Then $f$ is \emph{\(L\)-smooth} if
\begin{equation}
	\label{eq:gradsmooth}
\innp{\nabla f(y) - \nabla f(x)}{y - x} \leq L \norm{y - x}^{2} 	\quad \text{for all } x, y \in P,
\end{equation}
and similarly $f$ is \emph{$\mu$-strongly convex} if 
\begin{equation}
\innp{\nabla f(y) - \nabla f(x)}{y - x} \geq \mu \norm{y - x}^{2} 	\quad \text{for all } x, y \in P.
\end{equation}
\end{remark}

There is also the closely related and seemingly stronger property of \(L\)-Lipschitz continuous gradient \(\nabla f\), however in the case that $P$ is full-dimensional and $f$ is convex it is known to be equivalent to \(L\)-smoothness (see \citet[Theorem~2.1.5]{nesterov18} for the unbounded case, i.e., where \(P = \RR^{n}\) and \citet[Lemma 1.7]{CGFWSurvey2022} for $P$ being arbitrary convex domain). In particular, for twice differentiable convex functions $f$, we can also capture smoothness and strong convexity in terms of the Hessian via \(\norm{\nabla^{2} f} \leq L\) and via the largest eigenvalue of $\nabla^{2} f$ being upper bounded by $L \geq 0$ and the smallest eigenvalue being lower bounded by $\mu \geq 0$, respectively; the first inequality is useful for numerical estimation of $L$. 

In the following the domain $P$ will be a compact convex set and we assume that we have access to a so-called \emph{Linear Minimization Oracle (LMO)} for $P$, which upon being provided with a linear objective function $c$ returns a minimizer \(v = \argmin_{x\in P} \innp{c}{x}\) as formalized in Algorithm~\ref{ora:lmo}. Note that $v$ is not necessarily unique and without loss of generality we assume that $v$ is an extreme point of $P$; these extreme points are also often called \emph{atoms} in the context of Frank-Wolfe algorithms. For the compact convex set $P$ the \emph{diameter} $D$ of
$P$ is defined as $D \doteq \max_{x,y \in P} \norm{x - y}$.

Regarding the function $f$ we assume that we have access to gradients and function evaluations which is formalized as \emph{First-Order Oracle}\index{First-Order Oracle} (denoted as FOO), which given a point \(x \in P\), returns the function value \(f(x)\) and gradient \(\nabla f(x)\) at \(x\); see Algorithm~\ref{ora:FOO}. In the following we let $x^*$ be one of the optimal solutions to $\min_{x \in P} f(x)$ and define further $f^* \doteq f(x^*)$. Moreover, if not stated otherwise we consider $f\colon P \rightarrow \RR$. 

\begin{figure}[H]
	% \label{fig:algs}   % todo 
	\begin{minipage}[t]{0.49\textwidth}
		\centering
		\begin{oracle}[H]
			\label{ora:lmo}
			\caption{Linear Minimization Oracle over $P$ (LMO)}
			\begin{algorithmic}
				\REQUIRE Linear objective \(c\)
				\ENSURE \(v \in \argmin_{x\in P}\innp{c}{x}\)
			\end{algorithmic}
		\end{oracle}
	\end{minipage} 
	\hfill	
	\begin{minipage}[t]{0.49\textwidth}
		\centering
		\begin{oracle}[H]
			\caption{First-Order Oracle for \(f\) (FOO) \phantom{ over $P$ (LMO)}}
			\label{ora:FOO}
			\begin{algorithmic}
				\REQUIRE Point \(x \in P\)
				\ENSURE \(\nabla f(x)\) and \(f(x)\) \phantom{$\argmin_{x\in P}\innp{c}{x}$}
			\end{algorithmic}
		\end{oracle}
	\end{minipage}
	\caption{Access to function $f$ and feasible region $P$ is via two functions that we assume to have (oracle) access to.}
\end{figure}

\section{The Frank-Wolfe algorithm\pagehint{1}}
\label{sec:fw-algo}

We will now introduce the original variant of the Frank-Wolfe  (FW) algorithm due to \citet{fw56}, which is often also referred to as Conditional Gradients
\citep{polyak66cg}. Although many advanced variants with enhanced properties and improved convergence in specific problem configurations exist today, we will focus on the original version for clarity and to underscore the fundamental concepts. 

Suppose we are interested in minimizing a smooth and convex function $f$ over some compact convex feasible set $P$. A natural strategy would be to follow the negative of the gradient $\nabla f(x)$ at a given point $x$. However, how far can we go into that direction before we hit the boundary of the feasible region? Moreover, even if we would know how far we can go, i.e., we would potentially truncate steps to not leave the feasible region, even then the resulting algorithm might not be converging to an optimal solution. In fact, the arguably most well-known strategy, the \emph{projected gradient descent} method does not simply stop at the boundary but follows the negative of the gradient according to some step-size, disregarding the constraints, and then \emph{projects back} onto the feasible region. This last step can be very costly: if we do not have an efficient formulation or algorithm for the projection problem, solving this projection problem can be a (relatively expensive) optimization problem in itself. In contrast, the basic idea of the Frank-Wolfe algorithm is to \emph{not} follow the negative of the gradient but to follow an alternative direction of descent, which is well-enough aligned with the negative of the gradient, ensures enough primal progress, and for which we can easily ensure feasibility by means of computing convex combinations. This is done via the aforementioned Linear Minimization Oracle, with which we can optimize the negative of the gradient over the feasible region $P$ and then take the obtained vertex to form an alternative direction of descent. The overall process is outlined in Figure~\ref{fig:schematic_CG_step} and in Algorithm~\ref{alg:fw} we provide the Frank-Wolfe algorithm, which only requires access to \eqref{eq:mainProb} via the LMO (see Algorithm~\ref{ora:lmo}) to access the feasible region and via the FOO (see Algorithm~\ref{ora:FOO}) to access the function. 

\begin{algorithm}[h]
	\caption{Frank–Wolfe algorithm}
	\label{alg:fw}
	\begin{algorithmic}[1]
		\REQUIRE Initial atom $x_0\in P$,
		smooth and convex objective function \(f\)
		\ENSURE Iterates $x_1, \dotsc \in P$
		\FOR{$t=0$ \TO \dots}
		\STATE\label{fw_extract}
		$v_{t} \gets
		\argmin_{v \in P} \innp{\nabla f(x_{t})}{v}$
		\STATE \label{line:basic-step-size} \(\gamma_{t} \gets \frac{2}{2+t}\)
		\STATE $x_{t+1} \gets (1 - \gamma_t) x_{t} + \gamma_{t} v_{t}$
		\ENDFOR
	\end{algorithmic}
\end{algorithm}

As can be seen, assuming access to the two oracles, the actual implementation is very straight-forward: a simple computation of a convex combination, which ensures that we do not leave the feasible region. We made the deliberate choice in Line~\ref{line:basic-step-size} of Algorithm~\ref{alg:fw} to use the most basic step-size strategy $\gamma_{t} = \frac{2}{2+t}$, the so-called \emph{open loop} or \emph{agnostic} step-size, as this makes the algorithm \emph{parameter-independent}, i.e., not requiring any function parameters or parameter estimations. In the worst-case, this step-size is not dominated by more elaborate strategies (such as, e.g., line search or short steps), however in many important special cases there are better choices. As this is crucial we will discuss this a little more in-depth in Section~\ref{sec:adaptive} and will also provide a new variant of an adaptive step-size strategy.

Another important property is that the algorithm is \emph{affine invariant}, i.e., problem rescaling etc.~does not affect the algorithm's performance, compared to most other methods including PGD (notable exceptions exist, e.g., Newton's method). This makes the algorithm also very robust (especially with the open loop step-sizes) often offering superior numerical stability. 

Finally, we would like to mention that at iteration $t$ the iterate $x_t$ is a convex combination of at most $t+1$ extreme points (or atoms) of $P$. This will allow us later to obtain sparsity vs.~approximation trade-offs in Section~\ref{sec:convergence}.

\section{Properties\pagehint{8}}
\label{sec:properties}

We will now establish key properties of Algorithm~\ref{alg:fw}. We start with convergence properties and will then establish matching lower bounds as well as other properties. 

\subsection{Convergence}
\label{sec:convergence}

We will now prove the convergence of the Frank-Wolfe algorithm (Algorithm~\ref{alg:fw}). Convergence proofs for these methods typically use two key ingredients, which we will introduce in the following. 

\begin{lemma}[Primal gap, Dual gap, and Frank-Wolfe gap]
	\label{lem:gaps}
Let $f$ be a convex function and $P$ a compact convex set and consider \eqref{eq:mainProb}. For all $x \in P$ it holds:
\begin{equation}
	\tag{FW-gap}
	\label{eq:fwgap}
\underbrace{f(x) - f(x^*)}_{\text{primal gap at $x$}} \leq \underbrace{\innp{\nabla f(x)}{x - x^*}}_{\text{dual gap at $x$}} \leq \underbrace{\max_{v \in P} \innp{\nabla f(x)}{x - v}}_{\text{Frank-Wolfe gap at $x$}}.
\end{equation}
% The right-most quantity is the \emph{Frank-Wolfe gap at $x$}.
\begin{proof}
The first inequality follows from convexity and the second inequality follows from maximality.
\end{proof}
\end{lemma}

The Frank-Wolfe gap plays a crucial role in the theory of Frank-Wolfe methods as it provides an easily computable optimality certificate and suboptimality gap measure. An extreme point $v \in \argmax_{z \in P} \innp{\nabla f(x)}{x - z}$, is typically referred to as \emph{Frank-Wolfe vertex for $\nabla f(x)$}. The Frank-Wolfe gap also naturally appear in the first-order optimality condition for \eqref{eq:mainProb}, which states that $x^* \in P$ is optimal for \eqref{eq:mainProb} if and only if the Frank-Wolfe gap at $x^*$ is equal to $0$. Note that in the constrained case it does not necessarily hold that $\nabla f(x^*) = 0$.

\begin{lemma}[First-order Optimality Condition]
	\label{lem:foOptimality}
	Let $x^* \in P$. Then $x^*$ is an optimal solution to \eqref{eq:mainProb} if and
	only if 
	$$\innp{\nabla f(x^*)}{x^*-v} \leq 0$$ 
	for all $v \in P$. In particular, we have that the Frank-Wolfe gap
 $\max_{v \in P} \innp{\nabla f(x^*)}{x^*-v} = 0$.
\end{lemma}

The second property that is crucial is smoothness as it allows us to lower bound the primal progress we can derive from a step of the Frank-Wolfe algorithm.

\begin{lemma}[Primal progress from smoothness]
	\label{lem:primalProgessSmooth}
Let $f$ be an $L$-smooth function and let $x_{t+1} = (1 - \gamma_t) x_t + \gamma_t v_t$ with $x_t, v_t \in P$. Then we have
\begin{equation}
	\label{eq:smoothProgress}
	f(x_t) - f(x_{t+1}) \geq \gamma_t \innp{\nabla f(x_t)}{x_t - v_t} - \gamma_t^2 \frac{L}{2} \norm{x_t - v_t}^2.
\end{equation}
\begin{proof} The statement follows directly from the smoothness inequality \eqref{eq:baseSmooth}
\begin{align*}
	f(y) - f(x) \leq \innp{\nabla f(x)}{y-x}
	+ \frac{L}{2} \norm{y-x}^{2},
\end{align*}
choosing $x \leftarrow x_t$ and $y \leftarrow x_{t+1}$, plugging in the definition of $x_{t+1}$, and rearranging. This gives the desired inequality
\begin{equation*}
	f(x_t) - f(x_{t+1}) \geq \gamma_t \innp{\nabla f(x_t)}{x_t - v_t} - \gamma_t^2 \frac{L}{2} \norm{x_t - v_t}^2.
\end{equation*}
\end{proof}
\end{lemma}

With these two key ingredients (Lemma~\ref{lem:gaps} and Lemma~\ref{lem:primalProgessSmooth}) we can now establish the basic convergence rate of the Frank-Wolfe algorithm: 

\begin{theorem}[Primal convergence of the Frank-Wolfe algorithm]
\label{thm:fw-primal-conv}
Let $f$ be an $L$-smooth convex function and let $P$ be a compact convex set of diameter $D$. Consider the iterates of Algorithm~\ref{alg:fw}. Then the following holds:
\begin{equation*}
f(x_t) - f(x^*) \leq \frac{2LD^2}{t+2},
\end{equation*}
and hence for any accuracy $\varepsilon > 0$ we have $f(x_t) - f(x^*) \leq \varepsilon$ for all $t \geq \frac{2LD^2}{\varepsilon}$.

\begin{proof}
The convergence proof of the Frank-Wolfe algorithm follows an approach that is quite representative for convergence results in that area. The proof follows the template outlined in \citet{CGFWSurvey2022} and mimics closely the proof in \citet{jaggi13fw}.

Our starting point is the inequality from Lemma~\ref{lem:primalProgessSmooth}
\begin{align*}
	f(x_t) - f(x_{t+1}) & \geq \gamma_t \innp{\nabla f(x_t)}{x_t - v_t} - \gamma_t^2 \frac{L}{2} \norm{x_t - v_t}^2 \\ 
	& \geq \gamma_t (f(x_t) - f(x^*)) - \gamma_t^2 \frac{L}{2} \norm{x_t - v_t}^2. & \text{(Lemma~\ref{lem:gaps})}
\end{align*}
Subtracting $f(x^*)$ on both sides, bounding $\norm{x_t - v_t} \leq D$, and rearranging leads to
\begin{align}
	\label{eq:contractionFW}
	f(x_{t+1}) - f(x^*) \leq (1 - \gamma_t) (f(x_{t}) - f(x^*)) + \gamma_t^2 \frac{LD^2}{2}.
\end{align}
This contraction relates the primal gap at $x_{t+1}$ with the primal gap at $x_t$. We conclude the proof by induction. First observe that for $t = 0$ by \eqref{eq:contractionFW} it follows
$$
f(x_{1}) - f(x^*) \leq \frac{LD^2}{2} \leq \frac{2 LD^2}{2}.
$$
Now consider $t \geq 1$. We have 
\begin{align*}
	f(x_{t+1}) - f(x^*) & \leq (1 - \gamma_t) (f(x_{t}) - f(x^*)) + \gamma_t^2 \frac{LD^2}{2} \\
			& \leq \frac{t}{2+t} (f(x_{t}) - f(x^*)) + \frac{4}{(2+t)^2} \frac{LD^2}{2} & \text{(definition of $\gamma_t$)} \\
			& \leq \frac{t}{2+t} \frac{2LD^2}{2+t} + \frac{4}{(2+t)^2} \frac{LD^2}{2} & \text{(induction hypothesis)} \\
			& = \frac{2LD^2}{t+3} \left(\frac{(3+t)(1+t)}{(2+t)^2} \right) \leq \frac{2LD^2}{t+3}, & \text{($(3+t)(1+t) \leq (2+t)^2$)}
\end{align*}
which completes the proof.
\end{proof}
\end{theorem}

The theorem above provides a convergence guarantee for the primal gap. However, it relies on knowledge of the diameter $D$ and Lipschitz constant $L$ for estimating the number of required iterations to reach a certain target accuracy $\varepsilon$. We can also consider the Frank-Wolfe gap $\max_{v_t \in P} \innp{\nabla f(x_t)}{x_t - v_t}$, which upper bounds the primal gap $f(x_t) - f(x^*)$ via Lemma~\ref{lem:gaps}. While this gap is not monotonously decreasing (similar to the primal gap in the case of the open loop step-size) it is readily available in each iteration and hence can be used as a stopping criterion, i.e., we stop the algorithm when $\max_{v_\tau \in P} \innp{\nabla f(x_\tau)}{x_\tau - v_\tau} \leq \varepsilon$, not requiring a priori knowledge about $D$ and $L$. For the running minimum we can establish a convergence rate similar to that in Theorem~\ref{thm:fw-primal-conv}; see \citet{jaggi13fw}, see also \citet[Theorem 2.2 and Remark 2.3]{CGFWSurvey2022}. 

\begin{theorem}[Frank-Wolfe gap convergence of the Frank-Wolfe algorithm]
	\label{thm:dualConv}
Let $f$ be an $L$-smooth convex function and let $P$ be a compact convex set of diameter $D$. Consider the iterates of Algorithm~\ref{alg:fw}. 
Then the running minimum of the Frank–Wolfe gaps up to iteration $t$ satisfies:
\begin{equation*}
	\min_{0 \leq \tau \leq t} \max_{v_\tau \in P} \innp{\nabla f(x_\tau)}{x_\tau - v_\tau} \leq \frac{6.75\, L D^{2}}{t + 2}
\end{equation*}
\end{theorem}

Another important property of the Frank-Wolfe algorithm is that it maintains convex combinations of extreme points and in each iteration at most one new extreme point is added. This leads to a natural accuracy vs.~sparsity trade-off, where sparsity broadly refers to having convex combinations with a small number of vertices. This property is very useful and has been exploited repeatedly to prove mathematical results via applying the convergence guarantee of the Frank-Wolfe algorithm; we will see such an example further below in Section~\ref{sec:approxCara}

\subsection{A matching lower bound}
\label{sec:lowerbound}

In this section we will now provide a matching lower bound example due to \citet{lan2013complexity,jaggi13fw} that will require \(\Omega(L D^{2} / \varepsilon)\) LMO calls to achieve an accuracy of $\varepsilon$ for an $L$-smooth function $f$ and a feasible region of diameter $D$. This lower bound holds for \emph{any} algorithm that accesses the feasible region solely through an LMO and shows that in general the convergence rate of the Frank-Wolfe algorithm in Theorem~\ref{thm:fw-primal-conv} cannot be improved. We consider
\begin{equation*}
	\min_{x \in \conv{\{e_{1}, \dots, e_{n} \}}} \norm{x}^{2},
\end{equation*}
i.e., we minimize the standard quadratic
\(f(x) = \norm{x}^{2}\)
over the probability simplex
\(P \doteq \conv{\{e_{1}, \dots, e_{n}\}}\), where the $e_{i}$ denote the standard basis vectors in $\R^{n}$, i.e., we have $L=2$ and $D = \sqrt{2}$ and any other combination of values for $L$ and $D$ can be obtained via rescaling. As $f$ is strongly convex it has a unique optimal solution, which is easily seen to be $x^* = (\frac{1}{n}, \dots, \frac{1}{n})$ with optimal objective function value $f(x^*) = \frac{1}{n}$. Note that the optimal solution lies in the relative interior of $P$, one of the earliest cases in which improved convergence rates for Frank-Wolfe methods have been obtained \citep{gm86}. 

If we now run the Frank-Wolfe algorithm from any extreme point $x_0$ of $P$, then after $t < n$ iterations, we have made $t$ LMO calls, and hence have picked up at most $t+1$ of the $n$ standard basis vectors. This is the only information available to us about the feasible region and by convexity the only feasible points the algorithm can create are convex combinations $x_t$ of these picked up extreme points. Thus it holds
\begin{equation*}
	f(x_t) \geq \min_{\substack{x \in \conv{\mathcal{S}} \\
			\mathcal{S} \subseteq \{e_1, \dots, e_n\} \\
			\card{\mathcal{S}} \leq t + 1}}
	f(x) = \frac{1}{t + 1}.
\end{equation*}
Therefore the primal gap after $t$ iterations satisfies
\(f(x_{t}) - f(x^{*}) \geq 1 / (t+1) - 1/n\) and thus with the choice \(n \gg 1 / \varepsilon\) we need \(\Omega(1 / \varepsilon)\) LMO calls
to guarantee a primal gap of at most \(\varepsilon\). Finally, observe that this example also provides an inherent
sparsity~vs.~optimality tradeoff: if we aim for a solution with sparsity \(t\), then the primal gap can be as large as $f(x_t)-f(x^*) \geq 1 / (t+1) - 1/n$. 

\begin{figure}
	\begin{minipage}[t]{0.49\textwidth}
		\centering
		\includegraphics[scale=0.4]{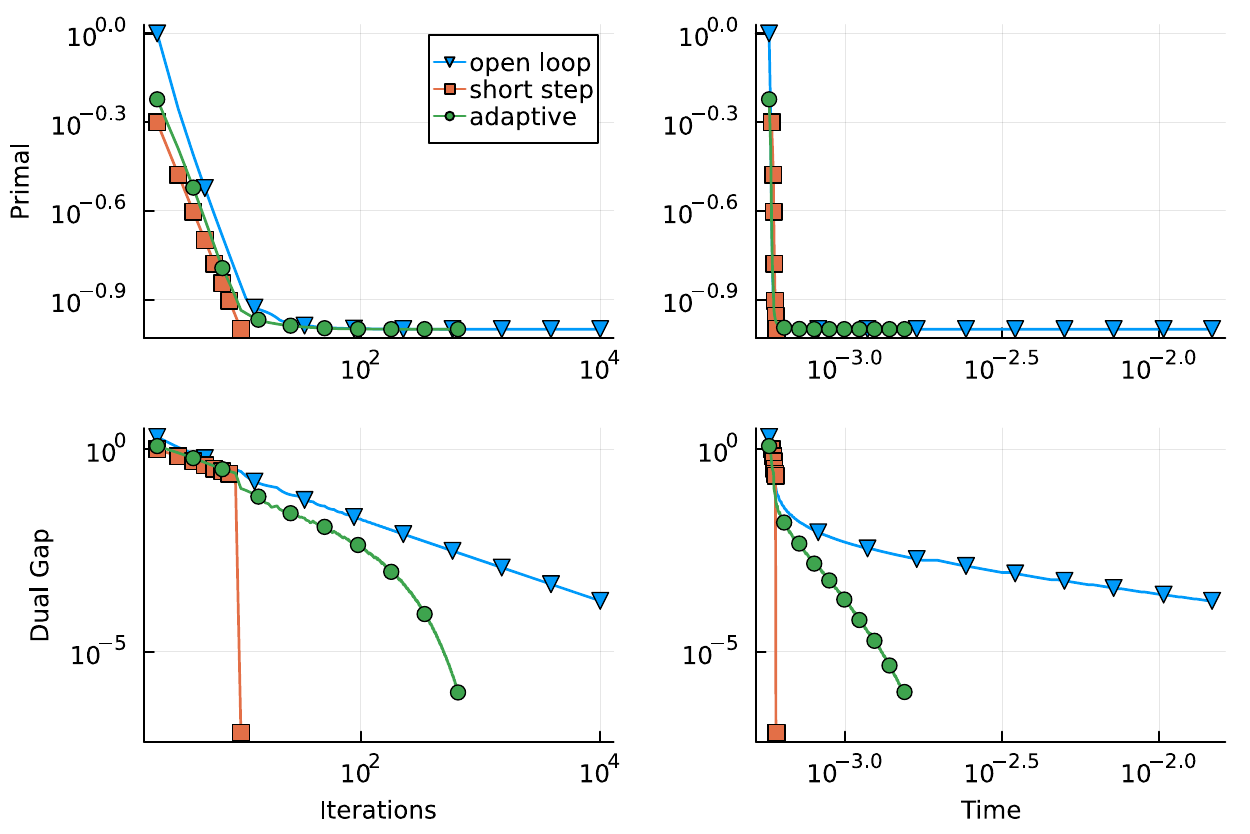}
		\caption{Minimizing $f(x) = \norm{x}^2$ over the probability simplex of dimension $n=10$ with an iteration limit of $k = 10^4$. It can be seen that once the iteration $t$ crosses the dimension threshold $n$ the short step strategy immediately recovers the optimal solution. 
		} 
		\label{fig:lb1}
	\end{minipage} 
	\hfill	
	\begin{minipage}[t]{0.49\textwidth}
		\centering
		\includegraphics[scale=0.4]{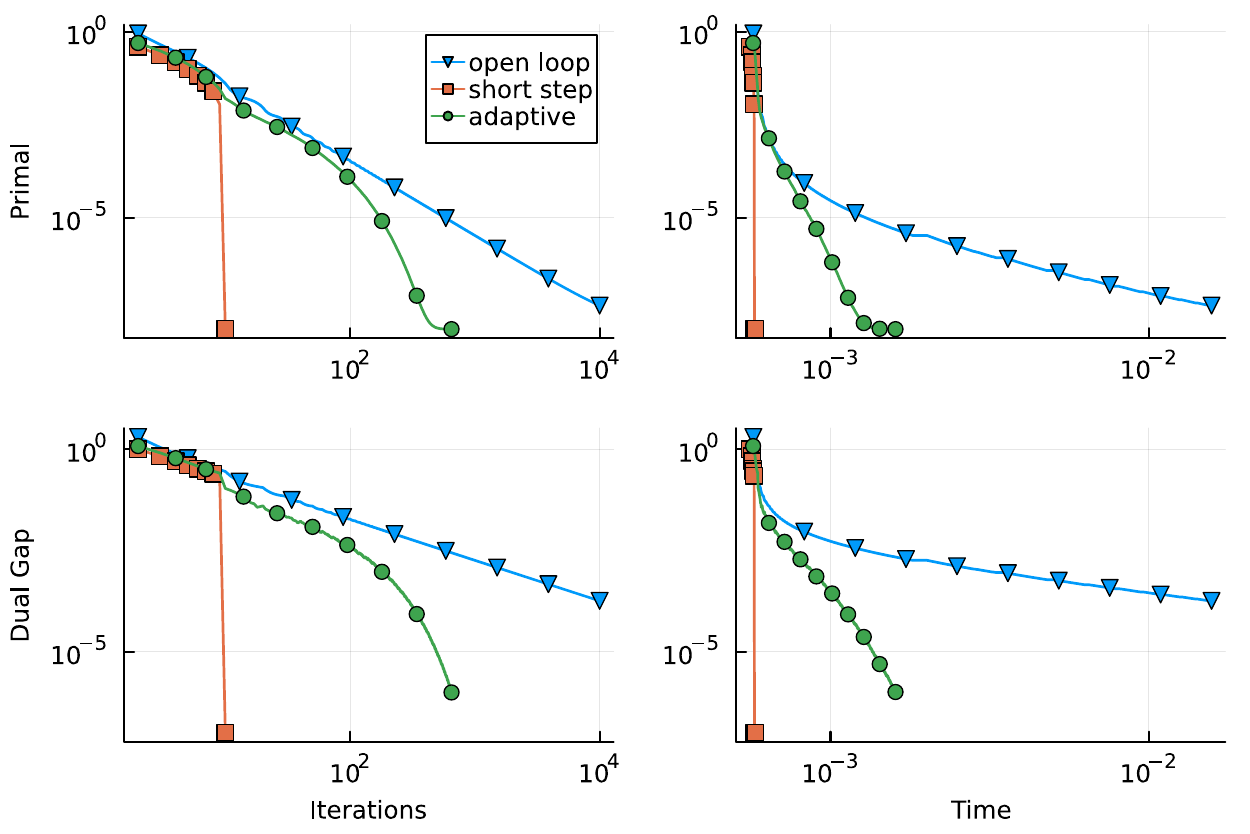}
		
		\caption{Same parameters as in Figure~\ref{fig:lb1}, however with the modified objective $f(x) = \norm{x - (\frac{1}{n}, \dots, \frac{1}{n})}^2$. Note, that $(\frac{1}{n}, \dots, \frac{1}{n})$ is the optimal solution for minimizing $f(x) = \norm{x}^2$ over the probability simplex. Dual convergence is identical while primal convergence differs.
		} 
		\label{fig:lb2}
	\end{minipage}
	
\end{figure}

However, several remarks are in order to put this example into perspective. First of all, the lower bound example only holds up to the dimension $n$ of the problem and that is for good reason. Once we pass the dimension threshold, the lower bound is not instructive any more and other step-size strategies might achieve linear rates for $t \geq n$, and in particular if the step-size is the short step rule (see also Section~\ref{sec:adaptive}) with exact smoothness $L$ we are optimal after exactly $t = n-1$ iterations; see Figures~\ref{fig:lb1} and~\ref{fig:lb2} for computational examples. Moreover, here we considered convergence rates independent of additional problem parameters. Introducing such parameters might provide more granular convergence rates under mild assumptions as shown, e.g., in \citet{garber2020sparseFW}. There is also a different lower bound of \(\Omega(1 / \varepsilon)\) by \citet{wolfe70} (see also \citet[Theorem 2.8]{CGFWSurvey2022}) that is based around the so-called zigzagging phenomenon of the Frank-Wolfe algorithm and that holds beyond the dimension threshold. However, it only holds for step-size strategies---grossly simplifying---that are at least as good as the short step strategy and interestingly the open loop step-size strategy is not subject to this lower bound. This is no coincidence, as there are cases \citep{WKP2022,WPP2023} where the open loop step-size can achieve a convergence rate of $\mathcal{O}(1/\varepsilon^2)$ for instances that satisfy the condition of the lower bound of \citet{wolfe70}. Finally, there is a universal lower bound \citep[Proposition 2.9]{CGFWSurvey2022} that matches the improved $\mathcal{O}(1/\varepsilon^2)$ rate for the open loop step-size: 

\begin{proposition}
	\label{prop:FW-fixed-lower}
	Let \(f\) be an $L$-smooth and convex function
	over a compact convex set \(P\).
	Then for $t \geq 1$, the iterates of the Frank-Wolfe algorithm (Algorithm~\ref{alg:fw}) with any step sizes \(\gamma_{\tau}\) satisfy
	\begin{align}
		\label{eq:FW-fixed-lower-value}
		f(x_{t}) - f(x^{*})
		&
		\geq
		\prod_{\tau=1}^{t-1} (1 - \gamma_{\tau})
		\cdot
		\innp{\nabla f(x^{*})}{x_{1} - x^{*}},
	\end{align}
	and in particular for the open loop step-size rule  \(\gamma_{\tau} = 2 / (\tau+2)\) we have
	\begin{equation}
		\label{eq:FW-fixed-lower-example}
		f(x_{t}) - f(x^{*})
		\geq
		\frac{2}{t (t+1)}
		\cdot
		\innp{\nabla f(x^{*})}{x_{1} - x^{*}}
		.
	\end{equation}
\end{proposition}

Finally, in actual computations these lower bounds are rarely an issue as instances often possess additional structure and adaptive step-size strategies (see Section~\ref{sec:adaptive}) provide excellent computational performance without requiring any knowledge of problem parameters. 

\subsection{Nonconvex objectives}
\label{sec:nonconvex}

The Frank-Wolfe algorithm can also be used to obtain locally optimal solutions if $f$ is nonconvex but smooth. In this case, $x \in P$ is \emph{locally optimal} (or \emph{first-order critical}) if and only if the Frank-Wolfe gap at $x$ is $0$, i.e., $\max_{v \in P} \innp{\nabla f(x)}{x - v} = 0$. We will present a simple argument to establish convergence to a locally optimal solution, however the argument can be improved as done in \citet{lj16nonconvex}, which was also the first to establish convergence for smooth nonconvex objectives. In particular, our argument will use a constant step-size \(\gamma_t = \gamma \doteq \frac{1}{\sqrt{T+1}}\) which has the advantage that it is parameter-free, but we need to decide on the number of iterations $T$ ahead of time and the convergence guarantee only holds for the last iteration $T$ in contrast to so-called \emph{anytime guarantees} that hold in each iteration $t = 0, \dots, T$. Nonetheless, the core of the argument is identical and more clearly isolated that way. 

\begin{theorem}[Convergence for nonconvex objectives]
	\label{thm:nonconvex} Let $f$ be an $L$-smooth but not necessarily convex function and $P$ be a compact convex set of diameter $D$. Let $T \in \N$, then the iterates of the Frank-Wolfe algorithm (Algorithm~\ref{alg:fw}) with the step-size \(\gamma_t = \gamma \doteq \frac{1}{\sqrt{T+1}}\) satisfy:
\begin{align*}
G_T \doteq \min_{0 \leq t \leq T} \max_{v_t \in P} \innp{\nabla f(x_t)}{x_t - v_t} \leq \frac{\max\{2h_0,
	LD^2\}}{\sqrt{T+1}},
\end{align*}
where $h_0 \doteq f(x_0) - f(x^*)$ is the primal gap at $x_0$.
\begin{proof}
	
Our starting point is the primal progress bound at iterate $x_t$ from Lemma~\ref{lem:primalProgessSmooth}
\begin{equation*}
	\label{eq:simpleSmoothNonconvex}
	f(x_t) - f(x_{t+1}) \geq \gamma \innp{\nabla f(x_t)}{x_t - v_t} - \gamma^2
	\frac{L}{2} \norm{x_t - v_t}^2.
\end{equation*}
Summing up the above along the iterations \(t = 0, \dots, T\) and
rearranging gives
\begin{align*}
	\gamma  \sum_{t = 0}^{T} \innp{\nabla f(x_t)}{x_t - v_t} & \leq f(x_0) -
	f(x_{T+1}) + \gamma^2 \sum_{t = 0}^{T} \frac{L}{2} \norm{x_t - v_t}^2
	\\
	& \leq f(x_0) -
	f(x^*) + \gamma^2 \sum_{t = 0}^{T} \frac{L D^2}{2} = h_0 + \gamma^2 (T+1) \frac{L D^2}{2}.
\end{align*}
We divide by \(\gamma (T+1)\) on both sides to arrive at our final
estimation
\begin{align}
	\label{eq:simpleNonconvexEst}
	G_{T} \leq \frac{1}{T+1} \sum_{t = 0}^{T} \innp{\nabla f(x_t)}{x_t - v_t} \leq \frac{h_0}{\gamma (T+1)} + \gamma \frac{L D^2}{2},
\end{align}
and for \(\gamma = \frac{1}{\sqrt{T+1}}\) this yields
\begin{align}
	\label{eq:24}
	G_{T} \leq \frac{1}{T+1} \sum_{t = 0}^{T} \innp{\nabla f(x_t)}{x_t -
		v_t} \leq \frac{2h_0 + LD^2}{2\sqrt{T+1}} \leq \frac{\max\{2h_0,
		LD^2\}}{\sqrt{T+1}},
\end{align}
which completes the proof. 
\end{proof}	
\end{theorem}

Note that $G_T$ can be observed throughout the algorithm's run and can be used as a stopping criterion. Moreover, the convergence rate of $\mathcal O(1/\sqrt{T})$ is optimal; see \citet{CGFWSurvey2022} for a discussion. If we have knowledge about \(h_0\), \(L\), and \(D\) then the above estimation can be slightly improved while maintaining a constant step size rule. We revisit \eqref{eq:simpleNonconvexEst} and optimize for
\(\gamma\), to obtain \(\gamma = \sqrt{\frac{2h_0}{LD^2(T+1)}}\) and hence:
\begin{align}
	\label{eq:25}
	G_{T} \leq \frac{1}{T+1} \sum_{t = 0}^{T} \innp{\nabla f(x_t)}{x_t -
		v_t} \leq \sqrt{\frac{2h_0 LD^2}{{T+1}}} \leq \frac{\max\{2h_0,
		LD^2\}}{\sqrt{T+1}}.
\end{align}
In the right most estimation the two bounds from \eqref{eq:24} and
\eqref{eq:25} are identical, which is due to the relatively weak estimation of the very last inequality. In fact, the
difference between \eqref{eq:24} und \eqref{eq:25} is that in the former we
have the arithmetic mean between \(2h_0\) and \(LD^2\) as bound, i.e., \(G_T
\leq \frac{2h_0 + LD^2}{2} \frac{1}{\sqrt{T+1}}\), whereas in \eqref{eq:25} we
have the geometric mean of the two terms, i.e., \(G_T \leq \sqrt{2h_0 LD^2}
\frac{1}{\sqrt{T+1}}\); by the AMGM inequality the latter is smaller than the
former. In both cases, we can also turn the guarantees into anytime guarantees (with minor changes in constants) by using the step-size rules $\gamma_t = 1/\sqrt{t+1}$ and \(\gamma_t = \sqrt{\frac{2h_0}{LD^2(t+1)}}\), respectively, and then using the bound $\sum_{t=0}^{T-1} \frac{1}{\sqrt{t+1}} \leq 2 \sqrt{T} - 1$. Then
telescoping works analogously to the above with minor adjustments. Finally, note that in all estimation we do not only provide a
guarantee for the running minimum of the Frank-Wolfe gap but their averages in
fact and the former is a consequence of the latter.

\subsection{Dual prices}

Another very useful property of the Frank-Wolfe algorithm (and also its more complex extensions) is that we readily obtain dual prices for active constraints, as long as the LMO provides dual prices. Similar to linear optimization, the dual price of a constraint captures the (local) rate of change of the objective if the constraint is relaxed. This is in particular useful in, e.g., portfolio optimization applications and energy problems, where marginal prices of constraints can guide decisions of real-world decision makers. Here we will only consider dual prices at the optimal solution $x^*$ and we will only cover the basic case without any degeneracy. However we can also derive dual prices for approximately optimal solutions and we refer the interested reader to \citet{BP2021} for an in-depth discussion. 

Suppose that the feasible region $P$ is actually a polytope of the form \(P = \{z : A z \leq b\}\) with $A \in \R^{m \times n}$ and $b \in \R^n$. Let $x \in P$ be arbitrary. By strong duality we have that $v \in P$ is a minimizer for the linear program 
\(\min_{\{z: Az \leq b\}} \innp{\nabla f(x)}{z}\) if and only if there exist dual prices \(0 \leq \lambda \in \R^m\), so that 
\begin{equation}
\tag{LP-duality}
\label{eq:grad-strong-duality}
	\nabla f(x) = - \lambda A \qquad \text{and} \qquad  \innp{\nabla f(x)}{v} = \min_{\{z: Az \leq b\}} \innp{\nabla f(x)}{z} = - \innp{\lambda}{b},
\end{equation}
i.e., the dual prices together with the constraints certify optimality. It is well known that the second equation can be equivalently replaced by a complementary slackness condition that states \(\innp{\lambda}{b - A v} = 0\); it can be readily seen that \eqref{eq:grad-strong-duality} implies \(\innp{\lambda}{b - A v} = 0\) by rearranging and the other direction follows similarly. Now consider a primal-dual pair $(v,\lambda)$ that satisfies \eqref{eq:grad-strong-duality}. By definition $v$ is also a Frank-Wolfe vertex for $\nabla f(x)$, so that we immediately obtain
\begin{equation*}
	\label{eq:4}
	\innp{\nabla f(x)}{x-v} = - \innp{\lambda A}{x-v} = \innp{\lambda}{b - A x},
\end{equation*}
i.e., the \emph{Frank-Wolfe gap at $x$ is equal to the complementarity gap for $x$ given $\lambda$}; if the latter would be $0$ then complementary slackness would hold or equivalently the Frank-Wolfe gap would be $0$ and $(x,\lambda)$ would be an optimal primal-dual pair. This can be now directly be related to $\min_{\{z: A z \leq b\}} f(z)$ via \emph{Slater's (strong duality) condition of optimality}: \(x\) is optimal for \(\min_{\{z: A z \leq b\}} f(z)\) if and
only if \(x\) is optimal for \(\min_{\{z: A z \leq b\}} \innp{\nabla f(x)}{ z}\). This implies that if $x$ is an optimal solution to \(\min_{\{z: A z \leq b\}} f(z)\) then $(x,\lambda)$ will also satisfy \eqref{eq:grad-strong-duality}. Hence for an optimal solution $x$, the dual prices $\lambda$ valid for $v$ are also valid for $x$.

Given that the LMO for polytopes is often realized via linear programming solvers that compute dual prices as by-product, we readily obtain dual prices \(\lambda\) for the optimal solution \(x^*\) via the Frank–Wolfe vertex \(v\) for \(\nabla f(x^*)\).

\subsection{Adaptive Step-sizes}
\label{sec:adaptive}

The primal progress of a Frank-Wolfe step is driven by the smoothness inequality. Suppose $f$ is $L$-smooth, then using the definition of the Frank-Wolfe step, i.e., $x_{t+1} = (1-\gamma_t) x_t + \gamma_t v_t$ and Lemma~\ref{lem:primalProgessSmooth} provides:
\begin{equation}
	\label{eq:smoothTemp}
	f(x_t) - f(x_{t+1}) \geq \gamma_t \innp{\nabla f(x_t)}{x_t - v_t} - \gamma_t^2 \frac{L}{2} \norm{x_t - v_t}^2.
\end{equation}
Now rather than plugging in the open loop step-size, we can view the right-hand side as an expression in one variable $\gamma_t$ and maximize the right-hand side. This leads to the optimal choice
\begin{equation}
	\label{eq:standSmoothProg}
	\gamma_t = \frac{\innp{\nabla f(x_t)}{x_t - v_t}}{L \norm{x_t - v_t}^2} \qquad \text{and} \qquad f(x_t) - f(x_{t+1}) \geq \frac{\innp{\nabla f(x_t)}{x_t - v_t}^2}{2L \norm{x_t - v_t}^2}.
\end{equation}
Technically we can only form convex combinations if $\gamma_t \in [0,1]$, so that we have to truncate $\gamma_t \coloneq \min\left\{\frac{\innp{\nabla f(x_t)}{x_t - v_t}}{L \norm{x_t - v_t}^2},1\right\}$; observe that $\gamma_t \geq 0$ holds always as the we have that the Frank-Wolfe gap $\innp{\nabla f(x_t)}{x_t - v_t} \geq 0$. This step-size choice is often referred to as \emph{short step step-size} and is the Frank-Wolfe equivalent to steepest descent. In the case that the truncation is active, it holds that $\innp{\nabla f(x_t)}{x_t - v_t} \geq L \norm{x_t - v_t}^2$ and together with \eqref{eq:smoothTemp} it follows that we are in a regime where we converge linearly with 
\begin{equation*}
f(x_t) - f(x_{t+1}) \geq \innp{\nabla f(x_t)}{x_t - v_t}/2 \geq (f(x_t) - f(x^*)) / 2, 
\end{equation*}
i.e., the primal progress is at least half of the Frank-Wolfe gap and hence at least half of the primal gap. 

The short step strategy avoids the overhead of line searches, however unfortunately it requires knowledge of the smoothness constant $L$ or at least reasonably tight upper bounds of such. This issue is what \citet{pedregosa2018step} addressed in a very nice paper by dynamically approximating $L$. Rather than performing a traditional line search on the function value, the approximation of $L$ leads only to a slightly slower convergence rate by a constant factor, has only small overhead, and does not suffer the additive resolution issue of traditional line searches, where one can only get as accurate as the line search $\varepsilon$. In particular, this adaptive strategy allows to adapt to the potentially better local smoothness of $f$, rather than relying on a worst-case estimate; see \citet{CGFWSurvey2022} for an in-depth discussion.

In a nutshell, what \citet{pedregosa2018step} do is perform a multiplicative search for $L$ until the smoothness inequality 
\begin{equation*}
	\tag{adaptive}
	\label{eq:adaptive}
	f(x_t) - f(x_{t+1}) \geq \gamma_t \innp{\nabla f(x_t)}{x_t - v_t} - \gamma_t^2 \frac{M}{2} \norm{x_t - v_t}^2.
\end{equation*}
holds for the approximation $M$ of $L$ with $\gamma_t = \min\left\{\frac{\innp{\nabla f(x_t)}{x_t - v_t}}{M \norm{x_t - v_t}^2},1\right\}$ being the short step. 

Unfortunately, checking \eqref{eq:adaptive} in practice can be numerically very challenging as we mix function evaluations, gradient evaluations, and quadratic norm terms. Rather we present a new variant of the adaptive step-size strategy, where we rely on a different test for accepting the estimation $M$ of $L$: 
\begin{equation}
	\tag{altAdaptive}
	\label{eq:altAdaptive}
	\innp{\nabla f(x_{t+1})}{x_t - v_t} \geq 0,
\end{equation}
where $x_{t+1} = (1-\gamma_t) x_t + \gamma_t v_t$ as before with $\gamma_t = \min\left\{\frac{\innp{\nabla f(x_t)}{x_t - v_t}}{M \norm{x_t - v_t}^2},1\right\}$ being the short step for the estimation $M$, i.e., we only test (inner products with) the gradient $\nabla f$ at different points. Moreover, this test might provide additional primal progress as we discuss below. This leads to the adaptive step-size strategy given in Algorithm~\ref{alg:AdaptiveStepSize}, which is numerically very stable, however requires gradient computations (rather than function evaluations).

\begin{figure}
	\begin{minipage}[t]{0.49\textwidth}
		\centering
		\includegraphics[scale=0.4]{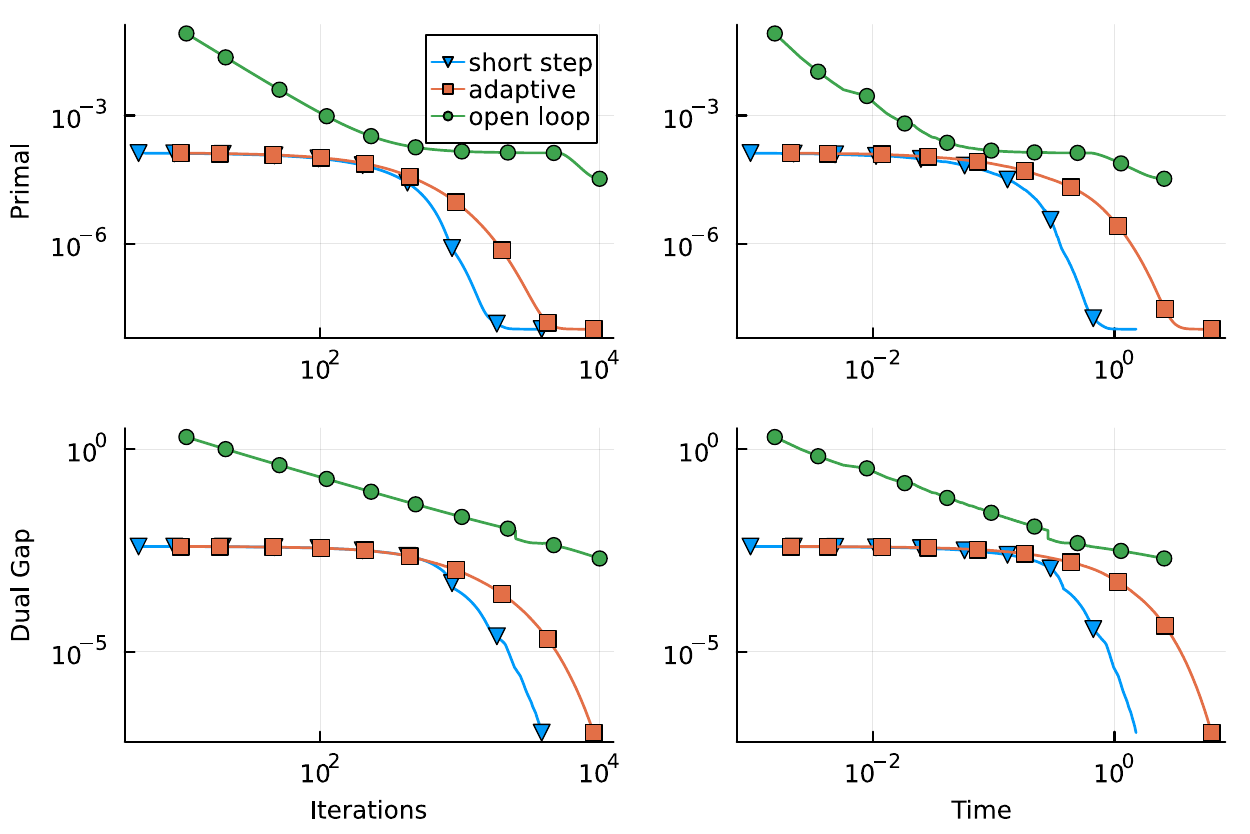}
		\caption{Convergence speed for a simple quadratic over a $K$-sparse polytope with three different step-size strategies. The basic open loop step-size $\gamma_t = \frac{2}{2+t}$, the short step rule, which requires a smoothness estimate $L$ (here we used the exact smoothness), and the adaptive step-size rule that dynamically approximates $L$. 			
		Plot is log-log so that the order of the convergence corresponds to different slopes of the trajectories.
		} 
		\label{fig:adaptive}
	\end{minipage} 
	\hfill	
	\begin{minipage}[t]{0.49\textwidth}
		\centering
		\includegraphics[scale=0.4]{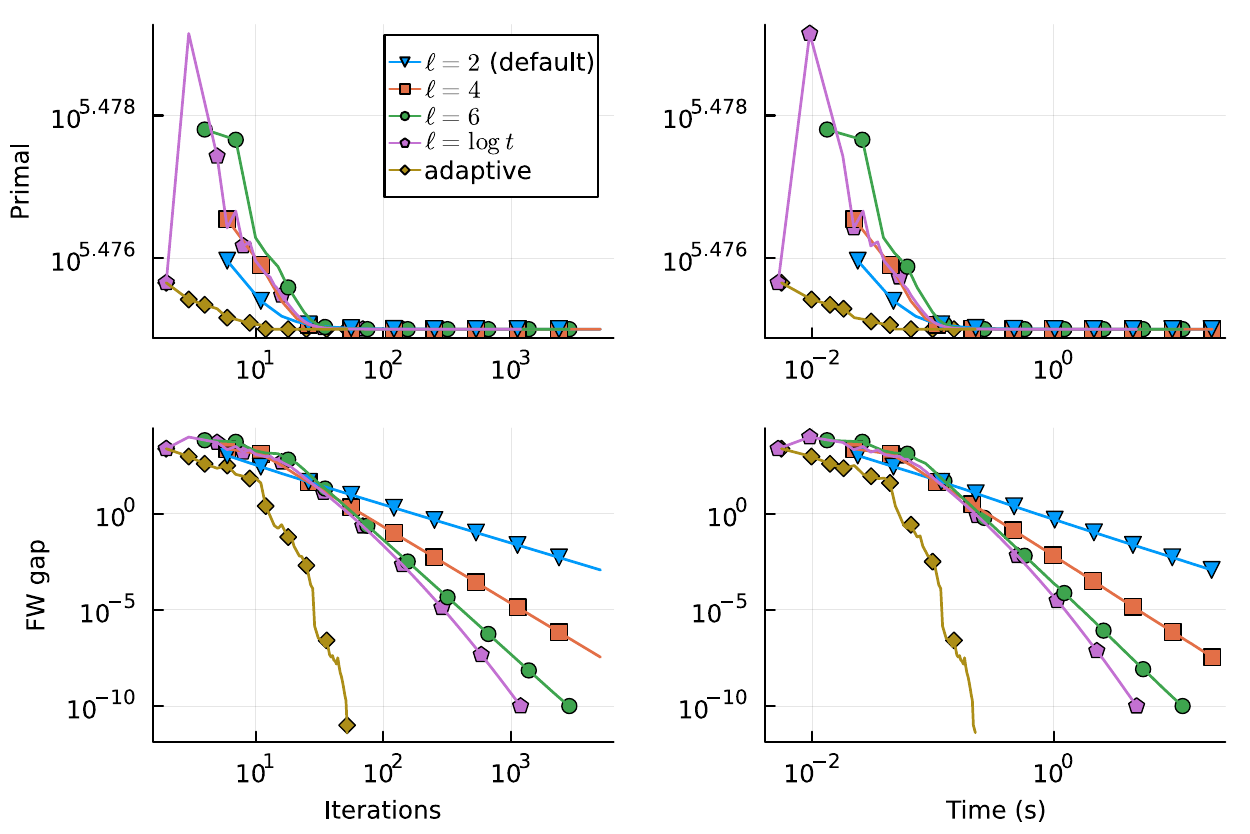}
		
		\caption{Convergence speed for a simple quadratic over a $K$-sparse polytope for open loop strategies of the form $\gamma_t = \frac{\ell}{\ell + t}$. We can see that (depending on the specifics of  the problem) larger $\ell$ achieve convergence rates of a higher order. For comparison also the adaptive step-size strategy has been included.
		Plot is log-log so that the order of the convergence corresponds to different slopes of the trajectories.
		} 
		\label{fig:open-loop}
	\end{minipage}
	
\end{figure}

We first show now that our condition \eqref{eq:altAdaptive} implies the same primal progress as \eqref{eq:adaptive} and then we will show that \eqref{eq:altAdaptive} holds for $L$ if $f$ is $L$-smooth. As such all results of \citet{pedregosa2018step} apply readily to the modified variant in Algorithm~\ref{alg:AdaptiveStepSize}. To demonstrate the convergence behavior of the various step-size strategies we ran a simple test problem with results presented in Figure~\ref{fig:adaptive}. We see that the adaptive strategy approximates the short step very well and both significantly outperform the open loop strategy. 

In the following we present the slightly more involved estimation based on a new progress bound from smoothness. For completeness we also include a significantly simplified estimation based on the regular smoothness bound in Appendix~\ref{sec:adaptive-easy}, however there we only guarantee approximation of the smoothness constant within a factor of $2$. We start with introducing another variant of the smoothness inequality. Note that all these inequalities are equivalent when considering all $x, y$, however we want to apply them for a \emph{specific} pair of points $x,y$ and then their transformations from one into another might not be sharp as demonstrated in the following remark: 

\begin{remark}[Point-wise smoothness estimations] Suppose that $f$ is $L$-smooth and convex and consider two points $x,y$. Suppose we want to derive \eqref{eq:baseSmooth} from the gradient-based variant in \eqref{eq:gradsmooth} using only the two points $x,y$. Then the naive way of doing so it:
\begin{align*}
f(y) - f(x) & \leq \innp{\nabla f(y)}{y - x} & \text{(convexity)} \\ 
			& \leq \innp{\nabla f(x)}{y - x} + L \norm{y-x}^2. & \text{(using \eqref{eq:gradsmooth})}
\end{align*}
Observe that this is \emph{almost} the desired inequality \eqref{eq:baseSmooth}, except for the smoothness constant $2L$ and not $L$.
\end{remark}

The following lemma provides a different smoothness inequality that allows for tighter estimations. It requires $f$ to be $L$-smooth and convex on a potentially slightly larger domain containing $P$. 

\begin{lemma}[Smoothness revisited]
\label{lem:smootRevisit}
Let $f$ be an $L$-smooth and convex function on the $D$-neighborhood of a compact convex set $P$, where $D$ is the diameter of $P$. Then for all $x,y \in P$ it holds:
\begin{equation}
\label{eq:strongSmooth}
\frac{\innp{\nabla f(y) - \nabla f(x)}{y-x}^2}{2L\norm{y-x}^2} \leq f(y) - f(x) - \innp{\nabla f(x)}{y-x}.
\end{equation}

\begin{proof}
As shown in \citet[Lemma 1.8]{CGFWSurvey2022}, if \(f\) is an \(L\)-smooth convex function on the \(D\)-neighborhood of a convex set \(P\),
then for any points \(x, y \in P\) it holds
\begin{equation}
		\label{eq:smooth-grad-diff}
\norm{\nabla f(y) - \nabla f(x)}^{2}
\leq
2L (f(y) - f(x) - \innp{\nabla f(x)}{y - x}).
\end{equation}
Next we lower bound the left-hand side as 
$$\frac{\innp{\nabla f(y) - \nabla f(x)}{y-x}^2}{\norm{y-x}^2} \leq \norm{\nabla f(y) - \nabla f(x)}^{2}.$$
Chaining these two inequalities together and rearranging gives the desired claim. 
\end{proof}
\end{lemma}

The proof above explicitly relies on the convexity of $f$ via \citet[Lemma 1.8]{CGFWSurvey2022}. With  Lemma~\ref{lem:smootRevisit} we can provide the following guarantee on the primal progress. 

\begin{lemma}[Primal progress from \eqref{eq:altAdaptive}]
Let $f$ be an $L$-smooth and convex function on the $D$-neighborhood of a compact convex set $P$, where $D$ is the diameter of $P$. Further, let $x_{t+1} = (1-\gamma_t) x_t + \gamma_t v_t$ with $\gamma_t = \min\left\{\frac{\innp{\nabla f(x_t)}{x_t - v_t}}{M \norm{x_t - v_t}^2},1\right\}$ for some $M$. If $\innp{\nabla f(x_{t+1})}{x_t - v_t} \geq 0$, then it holds:
\begin{equation*}
	f(x_t) - f(x_{t+1}) \geq 
	\begin{cases*}
		\frac{\innp{\nabla f(x_{t})}{x_{t} - v_{t}}^{2} +
			\innp{\nabla f(x_{t+1})}{x_{t} - v_{t}}^{2}}{2 \max\{L, M\}
			\norm{x_{t} - v_{t}}^{2}} & if \quad $\gamma_t \in [0, 1]$ \\
		\frac{\innp{\nabla f(x_{t})}{x_{t} - v_{t}}}{2} + \frac{
			\innp{\nabla f(x_{t+1})}{x_{t} - v_{t}}^{2}}{2 \innp{\nabla f(x_t)}{x_t - v_t}} \geq \frac{\innp{\nabla f(x_{t})}{x_{t} - v_{t}}}{2} & if \quad $\gamma_t = 1$ and $M \geq L$
	\end{cases*}.
\end{equation*}
\begin{proof}
	If $\gamma_t = 1$, then without loss of generality we can assume that $M = \frac{\innp{\nabla f(x_t)}{x_t - v_t}}{\norm{x_t - v_t}^2}$, as $M$ only occurs in the definition of $\gamma_t$ and $x_{t+1}$. 
	Our starting point is Equation~\eqref{eq:strongSmooth} with $x \leftarrow x_{t+1}$ and $y \leftarrow x_t$:
\begin{align*}
	f(x_{t}) - f(x_{t+1})
	& \geq
	\frac{\innp{\nabla f(x_{t}) - \nabla f(x_{t+1})}{x_{t} -
		x_{t+1}}^{2}}{2L \norm{x_{t} - x_{t+1}}^{2}}
	+ \innp{\nabla f(x_{t+1})}{x_{t} - x_{t+1}} \\
	& =
	\frac{\innp{\nabla f(x_{t}) - \nabla f(x_{t+1})}{x_{t} -
		v_{t}}^{2}}{2L \norm{x_{t} - v_{t}}^{2}}
	+ \frac{\innp{\nabla f(x_{t})}{x_{t} - v_{t}} \cdot
		\innp{\nabla f(x_{t+1})}{x_{t} - v_{t}}}{M \norm{x_{t} -
		v_{t}}^{2}} & \text{(definition of $x_{t+1}$ and $\gamma_t$)} \\
	& \geq
	\frac{\innp{\nabla f(x_{t}) - \nabla f(x_{t+1})}{x_{t} -
		v_{t}}^{2}}{2\max\{L, M\} \norm{x_{t} - v_{t}}^{2}}
	+ \frac{\innp{\nabla f(x_{t})}{x_{t} - v_{t}} \cdot
		\innp{\nabla f(x_{t+1})}{x_{t} - v_{t}}}{\max\{L, M\}
		\norm{x_{t} - v_{t}}^{2}} & \text{($\innp{\nabla f(x_{t+1})}{x_t - v_t} \geq 0$)} \\
	& =
	\frac{\innp{\nabla f(x_{t})}{x_{t} - v_{t}}^{2} +
		\innp{\nabla f(x_{t+1})}{x_{t} - v_{t}}^{2}}{2 \max\{L, M\}
		\norm{x_{t} - v_{t}}^{2}}.
\end{align*}
Now if $\gamma_t = 1$ and $M \geq L$, then the above simplifies to:
$$
f(x_{t}) - f(x_{t+1}) \geq \frac{\innp{\nabla f(x_{t})}{x_{t} - v_{t}}}{2} + \frac{
	\innp{\nabla f(x_{t+1})}{x_{t} - v_{t}}^{2}}{2 \innp{\nabla f(x_t)}{x_t - v_t}} \geq \frac{\innp{\nabla f(x_{t})}{x_{t} - v_{t}}}{2}.
$$
This finishes the proof.
\end{proof}
\end{lemma}

Before we continue a few remarks are in order. First of all, observe that
$$
f(x_{t}) - f(x_{t+1}) \geq \frac{\innp{\nabla f(x_{t})}{x_{t} - v_{t}}^{2} +
	\innp{\nabla f(x_{t+1})}{x_{t} - v_{t}}^{2}}{2 \max\{L, M\}
	\norm{x_{t} - v_{t}}^{2}},
$$ 
has an additional term compared to the standard smoothness estimation \eqref{eq:standSmoothProg} and $\innp{\nabla f(x_{t+1})}{x_{t} - v_{t}} = 0$ if and only if $x_{t+1}$ is identical to the line search solution. This is in particular the case if $f$ is a standard quadratic as then the line search solution is identical to the short step solution. Nonetheless, in the typical case this extra term provides additional primal progress. Taking the maximum in the denominator ensures that if $M < L$, then we recover the primal progress that one would have obtained with the estimation $M = L$. This seems counter-intuitive as usually using estimations $M < L$ would lead to overshooting and negative primal progress, however here we still require that \eqref{eq:altAdaptive} holds for $M$, which prevents exactly this as can be seen from the proof. In particular, disregarding adaptivity for a second, in the case where $L$ is known, then with the choice $M = L$, Lemma~\ref{lem:smootToAdaptive} provides a stronger primal progress bound compared to \eqref{eq:standSmoothProg} assuming that \eqref{eq:altAdaptive} holds for $L$ (which holds always as $f$ is $L$-smooth; see Lemma~\ref{lem:smootToAdaptive}):
$$
f(x_{t}) - f(x_{t+1}) \geq \frac{\innp{\nabla f(x_{t})}{x_{t} - v_{t}}^{2} +
	\innp{\nabla f(x_{t+1})}{x_{t} - v_{t}}^{2}}{2 L
	\norm{x_{t} - v_{t}}^{2}}.
$$
This improved primal progress bound might give rise to improved convergence rates in some regimes; see also \citet{teboulle2023elementary} for a similar analysis for the unconstrained case providing optimal constants for the convergence rates of gradient descent. Moreover, the discussion above also implies that if \eqref{eq:altAdaptive} holds it might provide more primal progress than the original test via \eqref{eq:adaptive} used in \citet{pedregosa2018step}.

To conclude, we will now show that \eqref{eq:altAdaptive} holds for $L$, whenever the function is $L$-smooth and $\gamma_t$ is the corresponding short step for $L$. This implies that both \eqref{eq:altAdaptive} and \eqref{eq:adaptive} hold for $L$ whenever $f$ is $L$-smooth with the added benefit of numerical stability and additional primal progress via \eqref{eq:altAdaptive}.

\begin{lemma}[Smoothness implies \eqref{eq:altAdaptive}] 
\label{lem:smootToAdaptive}	
Let $f$ be $L$-smooth. Further, let $x_{t+1} = (1-\gamma_t) x_t + \gamma_t v_t$ with $\gamma_t = \min\left\{\frac{\innp{\nabla f(x_t)}{x_t - v_t}}{L \norm{x_t - v_t}^2},1\right\}$. Then \eqref{eq:altAdaptive} holds, i.e.,
\begin{equation*}
\innp{\nabla f(x_{t+1})}{x_t - v_t} \geq 0.	
\end{equation*}
\begin{proof}
We use the alternative definition of smoothness using the gradients \eqref{eq:gradsmooth}, i.e., we have
$$\innp{\nabla f(y) - \nabla f(x)}{y - x} \leq L \norm{y - x}^{2} 	\quad \text{for all } x, y \in P.$$
Now plug in $x \leftarrow x_t$ and $y \leftarrow x_{t+1}$, so that we obtain
\begin{align*}
	\innp{\nabla f(x_{t+1}) - \nabla f(x_t)}{x_{t+1} - x_t} & \leq L \norm{x_{t+1} - x_t}^{2}
\end{align*}
and using the definition of $x_{t+1}$ it follows
\begin{align*}
\innp{\nabla f(x_{t+1}) - \nabla f(x_t)}{\gamma_t(v_t - x_t)} & \leq L \gamma_t^2 \norm{v_t - x_t}^{2}.
\end{align*}
Now, if $\gamma_t = 1$, then $\innp{\nabla f(x_t)}{x_t - v_t} \geq L \norm{x_t - v_t}^2$, so that 
$$
\innp{\nabla f(x_{t+1}) - \nabla f(x_t)}{v_t - x_t} \leq \innp{\nabla f(x_t)}{x_t - v_t},
$$
holds. Otherwise, if $0 < \gamma_{t} < 1$, then dividing by $\gamma_t$ and plugging in the definition of $\gamma_t$ yields
$$
\innp{\nabla f(x_{t+1}) - \nabla f(x_t)}{v_t - x_t} \leq \innp{\nabla f(x_t)}{x_t - v_t}.
$$
In both cases, rearranging gives the desired inequality 
$$
\innp{\nabla f(x_{t+1})}{x_t - v_t} \geq 0.
$$
Finally, in case $\gamma_{t} = 0$ we have $x_t = x_{t+1}$ and the assertion holds trivially.
\end{proof}
\end{lemma}

\begin{algorithm}[h]
	\caption[]{(modified) Adaptive step-size strategy}
	\label{alg:AdaptiveStepSize}
	\begin{algorithmic}[1]
		\REQUIRE Objective function \(f\),
		smoothness estimate \(\widetilde{L}\),
		feasible points \(x\), \(v\)
		with \(\innp{\nabla f(x)}{x - v} \geq 0\),
		progress parameters \(\eta \leq 1 < \tau\)
		\ENSURE Updated estimate \(\widetilde{L}^{*}\),
		step-size \(\gamma\)
		\STATE\(M \leftarrow \eta \widetilde{L}\)
		\LOOP
		\STATE
		\(\gamma \leftarrow
		\min \{\innp{\nabla f(x)}{x - v}
		\mathbin{/} (M \norm{x - v}^{2}), 1 \}\) \COMMENT{compute short step for estimation $M$}
		\IF{\(\innp{\nabla f(x + \gamma (v - x))}{x - v} \geq 0\)} 
		\label{SufficientDecrease}
		\STATE \(\widetilde{L}^{*} \leftarrow M\)
		\RETURN \(\widetilde{L}^{*}\), \(\gamma\)
		\ENDIF
		\STATE\(M \leftarrow \tau M\)
		\ENDLOOP
	\end{algorithmic}
\end{algorithm}

\begin{remark}[Faster open loop convergence]
The adaptive step-size strategy from above uses feedback from the function and as such is not of the open loop type. In many applications such adaptive strategies are the strategies of choice as the function feedback is relatively minimal but convergence speed is superior (in most but not all cases as mentioned in Section~\ref{sec:lowerbound}). For many important cases we can also obtain improved convergence with rates of higher order for open loop step-sizes by using the modified step-size rule $\gamma_t = \frac{\ell}{t + \ell}$ with $\ell \in \N_{\geq 1}$; see \citet{WKP2022,WPP2023} for details. This is quite surprising as we only change the shift $\ell$ and not the order of $t$ in the denominator of $\gamma_t$. In fact, note that the order of $t$ in the denominator cannot be changed significantly as we need that $\sum_t \gamma_t = \infty$ and that $\sum_t \gamma_t^2$ converges for the step-size strategy to work; see \citet{CGFWSurvey2022}. If the corresponding $\ell$ cannot be set in practice, and if it has to be an open loop strategy, $\gamma_t = \frac{2 + \log (t+1)}{t + 2 + \log (t+1)}$ works very well; we can use $t$ or $t+1$ in the log depending on whether the first iteration is $t = 0$ or $t=1$. This corresponds essentially to a strategy where $\ell$ is gradually increased and it provides accelerated convergence rates when those exist while maintaining the same worst-case convergence rates as the basic strategy $\gamma_t = \frac{2}{t + 2}$ \citep{WPP2023B}. For a sample computation, see Figure~\ref{fig:open-loop}. 
\end{remark}

\section{Two applications\pagehint{8}}
\label{sec:applications}

In the following we will present two applications of the Frank-Wolfe algorithm. Both examples use very simple quadratic objectives of the form $f(x) = \norm{x - p}^2$ for some $p$ for the sake of exposition; for more complex examples see also \citet{CGFWSurvey2022}. 

\subsection{The Approximate Carathéodory Problem\pagehint{2}}
\label{sec:approxCara}

Our first example, is the \emph{Approximate Carathéodory Problem}. For this example, the Frank-Wolfe algorithm does not only present a practical means to solve the problem but in fact, its convergence guarantees itself provide a proof of the theorem and optimal bounds for a wide variety of regimes. Here we will confine ourselves to the $2$-norm case not assuming any additional properties, however many more involved cases are possible as studied in \citet{cp2019approxCara}. 

Given a compact convex set $P \subseteq \RR^n$, recall that Carath\'eodory's theorem states that any $x^* \in P$ can be written as a convex combination of no more than $n+1$ extreme points of $P$, i.e., $x^* = \sum_{1 \leq i \leq n+1} \lambda_i v_i$ with $\lambda \geq 0$, $\sum_i \lambda_i = 1$, and $v_i$ extreme points of $P$ with $1 \leq i \leq n+1$. In the context of Carath\'eodory's theorem, the \emph{cardinality} of a point $x^* \in P$, refers to the minimum number of required extreme points to express $x^*$ as a convex combination of those. If $x^*$ is of low cardinality it is often also referred to as \emph{sparse}. Every specific convex combination that expresses $x^*$ provides an upper bound on the cardinality of $x^*$.
The approximate variant of Carath\'eodory's problem asks: given $x^* \in P$, what is required cardinality of an $x \in P$ to approximate $x^*$ within an error of no more than $\varepsilon > 0$ (in a given norm)? Put differently, we are looking for $x \in P$ with $\norm{x - x^*} \leq \varepsilon$ of low cardinality. The approximate Carath\'eodory theorem states:

\begin{theorem}[Approximate Carath\'eodory Theorem]
	\label{thm:appxCara}
Let $p\geq2$ and $P$ be a compact convex set. For every $x^*\in P$,
there exists $x \in P$ with
cardinality of no more than $\mathcal{O}(pD_p^2/\epsilon^2)$ satisfying
$\norm{x-x^*}_p \leq \epsilon$, where $D_p = \max_{v,w \in P}\norm{w-v}_p$ is the $p$-norm diameter of $P$.
\end{theorem}

Note that the bounds of Theorem~\ref{thm:appxCara} are essentially tight in many cases \citep{mirrokni2017tight}.
In the following, we briefly discuss the case $p = 2$ without any additional assumptions. Suppose we have given a point $x^* \in P$ we can consider the objective 
$$f(x) \doteq \norm{x - x^*}^2.$$ 
Further, let $\varepsilon > 0$ be the approximation guarantee. Assuming, we have access to an LMO for $P$, we can now minimize the function $f(x)$ over $P$  via the Frank-Wolfe algorithm (Algorithm~\ref{alg:fw}). In order to achieve $\norm{x - x^*} \leq \varepsilon$ we have to run the Frank-Wolfe algorithm until $f(x_t) = f(x_t) - f(x^*) \leq \varepsilon^2$, which by Theorem~\ref{thm:fw-primal-conv} takes no more than $\mathcal{O}(2 D^2/\epsilon^2)$ iterations, where $D$ is the $\ell_2$-diameter of $P$. Moreover, in each iteration the algorithm is picking up at most one extreme point as discussed in Section~\ref{sec:convergence}. This establishes the guarantee for case $p = 2$ in Theorem~\ref{thm:appxCara}. Here we applied the basic convergence guarantee from Theorem~\ref{thm:fw-primal-conv}. However, for the Frank-Wolfe algorithm many more convergence guarantees are known, depending on properties of the feasible domain and position of the point $x^*$ that we want to approximate with a sparse convex combination. These improved convergence rates immediately translate into improved approximation guarantees for the approximate Carathéodory problem and we state some of these guarantees in Table~\ref{tab:lit}. 

\begin{table}[h]
	\centering{
		\footnotesize
		\begin{tabular}{lll}
			\toprule
			\textbf{$\ell_p$-norm}&\textbf{Assumption}&{\textbf{Cardinality bound}}\\
			\midrule
		$p\in\left[2,+\infty\right[$&--&$\displaystyle\mathcal{O}\!\left(\frac{pD_p^2}{\epsilon^2}\right)$ or $\displaystyle\mathcal{O}\!\left(\frac{p(D_*^2+D_0^2)}{\epsilon^2}\right)$ 
		\\
			&& \hfill ($D_p, D_0, D_*$ diameters; see \citet{cp2019approxCara})\\
			\cmidrule(lr){2-3}
			& $x^*\in\operatorname{ri}(P) $&$\displaystyle\mathcal{O}\!\left(p\left(\frac{D_p}{r_p}\right)^2\ln\!\left(\frac{1}{\epsilon}\right)\right)$ 
			\\
			&&\hfill ($\operatorname{ri}(P)$ relative interior, $r_p$ radius so that $B^p_{r_p}(x^*) \cap \operatorname{aff}(P) \subseteq P$) \\
			\cmidrule(lr){2-3}
			&$\alpha_p$-strongly convex $P$  &$\displaystyle\mathcal{O}\!\left(\frac{\sqrt{p}D_p+p/\alpha_p}{\epsilon}\right)$ \\ 
			\cmidrule(lr){2-3}
			&$(\alpha_p,q_p)$-uniformly convex $P$, $q_p\in\left[2,+\infty\right[$ &$\displaystyle\mathcal{O}\!\left(\frac{(pD_p^2)^{(q_p-1)/q_p}+p/\alpha_p^{2/q_p}}{\epsilon^{2(q_p-1)/q_p}}\right)$ \\ 

			\midrule
			$p\in\left]1,2\right[$&--&$\displaystyle\mathcal{O}\!\left(\frac{n^{(2-p)/p}D_2^2}{\epsilon^2}\right)$ 
			\\
			\midrule
			$p=1$&--&$\displaystyle\mathcal{O}\!\left(\frac{nD_2^2}{\epsilon^2}\right)$ \\ % &--\\
			&&\hfill ($n$ ambient dimension of $P$, $D_2$ is $\ell_2$-diameter) \\ % 
			\midrule
			$p=+\infty$&--&$\displaystyle\mathcal{O}\!\left(\frac{D_2^2}{\epsilon^2}\right)$ \\ % 
			&&\hfill ($D_2$ is $\ell_2$-diameter) \\ % 
			\bottomrule
		\end{tabular}
	\caption{\label{tab:lit} Cardinality bounds to achieve $\epsilon$-approximation for the approximate Carath\'eodory problem with respect to the $\ell_p$-norm; see \citet{cp2019approxCara} for full table. 
Recall that $P$ is \emph{$(\alpha, q)$-uniformly convex}
if for any $x, y \in P$, $\gamma \in [0,1]$, and $z \in \R^n$ with $\norm{z} \leq 1$ we have $y + \gamma (x - y) + \gamma (1 - \gamma) \cdot \alpha \norm{x - y}^{q} z \in P$, where \(\alpha\) and \(q\) are positive. An \((\alpha/2, 2)\)-uniformly convex set is called \emph{\(\alpha\)-strongly convex}.
}
	}
\end{table}

\begin{figure}
	\begin{minipage}[t]{0.49\textwidth}
		\centering
		\includegraphics[scale=0.4]{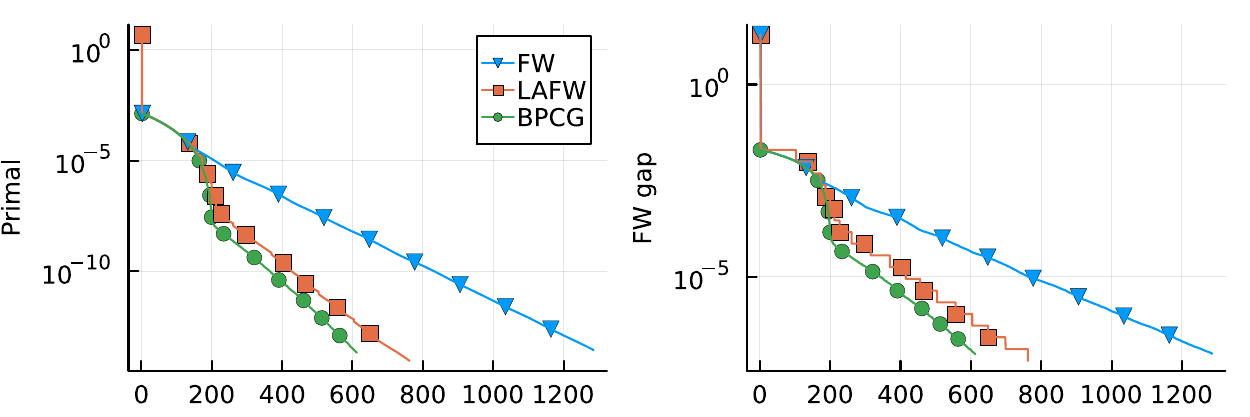}
		\caption{Cardinality vs.~approximation error in $\ell_2$-norm over a polytope of dimension $n = 1000$ for the Frank-Wolfe algorithm and two more advanced variants \emph{Lazy Away-step Frank-Wolfe} \citep{BPZ2017jour} and \emph{Blended Pairwise Conditional Gradients} \citep{TTP2021}. All algorithms use the short step step-size.  
		} 
		\label{fig:cara1}
	\end{minipage} 
	\hfill	
	\begin{minipage}[t]{0.49\textwidth}
		\centering
		\includegraphics[scale=0.4]{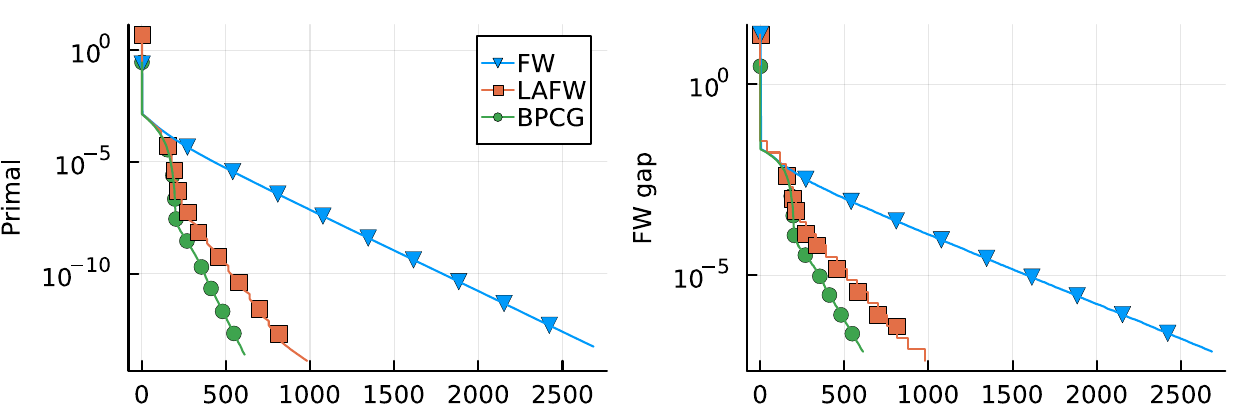}
		\caption{Same setup as in Figure~\ref{fig:cara1}, however the step-size strategy is the adaptive strategy from Section~\ref{sec:adaptive}. As we can see the Frank-Wolfe algorithm is quite sensitive to the strategy while more advanced variants, due to their design of only adding new vertices when not enough progress can be made otherwise, are not.
		} 
		\label{fig:cara2}
	\end{minipage}
	
\end{figure}

Apart from establishing theoretical results by means of the Frank-Wolfe algorithm's convergence guarantees it can be easily used in actual computations, see e.g., Figures~\ref{fig:cara1} and~\ref{fig:cara2} for an example. Observe that in the particular case of $f(x)$ here, we can also directly observe the primal gap and hence use it as a stopping criterion. The Frank-Wolfe approach to the approximate Carathéodory Problem has been also recently used in Quantum Mechanics to establish new Bell inequalities and local models, as well as improve the Grothendieck constant $K_G(3)$ of order $3$ (see \citet{DIBKGP2023,DVP2023} and the references contained therein). This approach is also very useful in the context of the \emph{coreset problem}, which asks for a subset of data points of a large data set that maintains approximately the same statistical properties (see \citet[Section 5.2.5]{CGFWSurvey2022}).  

\subsection{Separating hyperplanes\pagehint{3}}
\label{sec:sepa}

In Section~\ref{sec:approxCara} we have used the Frank-Wolfe algorithm for obtaining a convex decomposition of a point $x^* \in P$, i.e., we certified membership in $P$. We can also use the Frank-Wolfe algorithm with the same objective $f(x) = \norm{x - \tilde x}^2$ to obtain separating hyperplanes for points $\tilde x \not\in P$, i.e., we certify non-membership. This has been successfully applied in \citet{DIBKGP2023,DVP2023} to certify that the correlations of certain quantum states exhibit non-locality, i.e., are truly quantum by separating them from the local polytope, the polytope of all classical correlations. Moreover, it has also been used in \citet{TSPP2023} to derive separating hyperplanes from enumeration oracles. 

In the following we provide the most naive way of computing separating hyperplanes. An improved strategy has been presented in \citet{TSPP2023}, which derives a new algorithmic characterization of non-membership, that requires fewer iterations of the Frank-Wolfe algorithm compared to our naive strategy here. It is also interesting to note that from a complexity-theoretic perspective, what the Frank-Wolfe algorithm does is to turn an LMO for $P$ into a separation oracle for $P$ via optimizing the objective $\norm{x - \tilde x}^2$. 

Given $\tilde x \not\in P$, we consider the optimization problem
\begin{equation}
\label{eq:sep}
\tag{Sep}
\min_{x \in P} \norm{x - \tilde x}^2,
\end{equation}
with $f(x) \doteq \norm{x - \tilde x}^2$. Using Lemma~\ref{lem:foOptimality} we can immediately obtain a separating hyperplane from an optimal solution $x^* \in P$ to \eqref{eq:sep}:
\begin{equation}
	\label{eq:sepaHyper}
	\tag{sepHyperplane}
	\innp{\nabla f(x^*)}{x^*} \leq \innp{\nabla f(x^*)}{v},
\end{equation}
which holds for all $v \in P$. Moreover, as $\tilde{x} \not \in P$, we have by convexity $\innp{\nabla f(x^*)}{x^* - \tilde{x}} \geq f(x^*) - f(\tilde{x}) = f(x^*) > 0$. and hence \eqref{eq:sepaHyper} is violated by $\tilde{x}$, i.e., $\innp{\nabla f(x^*)}{x^*} > \innp{\nabla f(x^*)}{\tilde{x}}$. This argument provides the desired separating hyperplane mathematically, but numerically it is problematic as we usually solve Problem~\eqref{eq:sep} only up to some accuracy $\varepsilon > 0$, typically using the Frank-Wolfe gap $\max_{v \in P} \innp{\nabla f(x_t)}{x_t-v} \leq \varepsilon$ as stopping criterion. When the algorithm stops we similarly obtain 
\begin{equation}
	\label{eq:validIneq}
	\tag{validIneq}
	\innp{\nabla f(x_t)}{x_t} - \varepsilon \leq \innp{\nabla f(x_t)}{x} \quad \text{ which simplifies to }\quad \min_{v \in P} \innp{\nabla f(x_t)}{v} \leq \innp{\nabla f(x_t)}{x},
\end{equation}
which is valid for all $x \in P$. However this inequality does not necessarily separate $\tilde x$ from $P$. A sufficient condition for separation is that $\tilde x$ is $\sqrt{\varepsilon}$-far from $P$ so that we have $\norm{x^* - \tilde{x}} > \sqrt{\varepsilon}$. We then can use the same convexity argument as before: 
\begin{align*}
	& \innp{\nabla f(x_t)}{x_t - \tilde{x}} - \varepsilon & \text{(stopping criterion)} \\
	\geq &\; f(x_t) - f(\tilde{x}) - \varepsilon & \text{(convexity)} \\
	\geq &\; f(x^*) - \varepsilon > 0. & \text{($\norm{x^* - \tilde{x}} > \sqrt{\varepsilon}$)}
\end{align*}

Now turning this argument around, if $\tilde x \not \in P$ is $\varepsilon$-far from $P$, we need to run the Frank-Wolfe algorithm until the Frank-Wolfe gap satisfies $\max_{v \in P} \innp{\nabla f(x_t)}{x_t-v} \leq \varepsilon^2$. Combining this with Theorem~\ref{thm:dualConv} we can estimate
\begin{align*}
	\frac{6.75 LD^2}{t+2} \leq \varepsilon^2,
\end{align*}
with $L = 2$. Thus we have found a separating hyperplane for $\tilde x$ whenever $t \geq 13.5 D^2 / \varepsilon^2$. 

\begin{figure}[h]
	\begin{minipage}[b]{0.3\textwidth}
		\centering
		\includegraphics[scale=1]{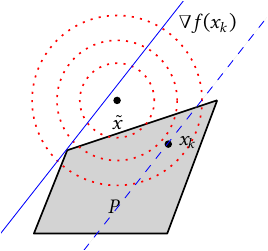}
	\end{minipage} 
	\hfill	
	\begin{minipage}[b]{0.3\textwidth}
		\centering
		\includegraphics[scale=1]{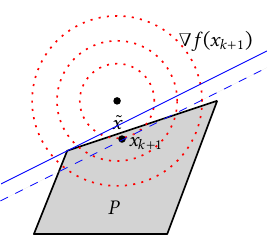}
	\end{minipage}
	\hfill
	\begin{minipage}[b]{0.3\textwidth}
		\centering
		\includegraphics[scale=1]{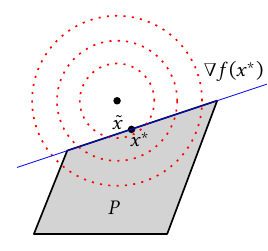}
		
	\end{minipage}
		\caption{\label{fig:sepa} 
			Each iterate $x_k$ induces a valid inequality of the form $\min_{v \in P} \innp{\nabla f(x_t)}{v} \leq \innp{\nabla f(x_t)}{x}$ in blue (dashed blue line is $\nabla f(x_t)$ at $x_t$), which may or may not separate $\tilde x$; see middle and left respectively. At $x^*$ the induced inequality $\min_{v \in P} \innp{\nabla f(x^*)}{v} = \innp{\nabla f(x^*)}{x^*}\leq \innp{\nabla f(x^*)}{x}$ is guaranteed to separate $\tilde x \not\in P$ and often (but not always) induces a facet of $P$.
		}	
\end{figure}

In practice however we usually do not know $D$ and we also do not know how far $\tilde x$ is from $P$. Nonetheless, we can simply test in each iteration $t$ whether $\tilde x$ violates \eqref{eq:validIneq}, i.e., 
$$
\nabla f(x_t) \text{ separates } \tilde x \quad \Leftrightarrow \quad  \min_{v \in P}  \innp{\nabla f(x_t)}{v} > \innp{\nabla f(x_t)}{\tilde x}.
$$
and simply stop then and are guaranteed this will take no more than $\mathcal O(D^2/\varepsilon^2)$ iterations. The process is illustrated in Figure~\ref{fig:sepa}. A similar approach, basically combining Sections~\ref{sec:approxCara} and~\ref{sec:sepa} can also be used to compute the intersection of two compact convex sets or certify their disjointness by means of a separating hyperplane (assuming LMO access to each) as shown in \citet{BPW2022}.

\section{Computational codes}
\label{sec:computations}

For actual computations we have developed the \href{https://github.com/ZIB-IOL/FrankWolfe.jl}{\texttt{FrankWolfe.jl}} Julia package, which implements many Frank-Wolfe variants and is highly customizable. Moreover, we have also developed a mixed-integer extension \href{https://github.com/ZIB-IOL/Boscia.jl}{\texttt{Boscia.jl}} that allows for some of the variables taking discrete values. 

\section*{Acknowledgments}
The author would like to thank Gábor Braun for pointing out the alternative smoothness inequality used in Section~\ref{sec:adaptive}, which gave rise to tighter bounds. This research was partially supported by the DFG Cluster of Excellence MATH+ (EXC-2046/1, project id 390685689) funded by the Deutsche Forschungsgemeinschaft (DFG).

\bibliographystyle{apalike}
\bibliography{bibs/bibliography,bibs/pubPokutta,bibs/cgfw_bibliography,bibs/web}

\clearpage
\appendix

\section{Adaptive Step-sizes: Simpler estimation}
\label{sec:adaptive-easy}

In this section we will present a simplified estimation of Section~\ref{sec:adaptive} for adaptive step-sizes albeit at the cost being only able to approximate the smoothness of $f$ within a factor of $2$. 
The basic setup is identical to the one before, however we use a different test for accepting the estimation $M$ of $L$: 
\begin{equation}
	\tag{altAdaptive-simple}
	\label{eq:altAdaptive-simple}
	\innp{\nabla f(x_{t+1})}{x_t - v_t} \geq \frac{1}{2} \innp{\nabla f(x_{t})}{x_t - v_t},
\end{equation}
where $x_{t+1} = (1-\gamma_t) x_t + \gamma_t v_t$ as before with $\gamma_t = \min\left\{\frac{\innp{\nabla f(x_t)}{x_t - v_t}}{M \norm{x_t - v_t}^2},1\right\}$ being the short step for the estimation $M$ and the corresponding algorithm becomes Algorithm~\ref{alg:AdaptiveStepSize-simple} in this case. We proceed similarly as before: we first show that condition \eqref{eq:altAdaptive-simple} implies primal progress and then we will show that \eqref{eq:altAdaptive-simple} holds for $2L$ if $f$ is $L$-smooth; this is were \eqref{eq:altAdaptive-simple} is weaker than \eqref{eq:altAdaptive}.

\begin{lemma}[Primal progress from \eqref{eq:altAdaptive-simple}]
	Let $x_{t+1} = (1-\gamma_t) x_t + \gamma_t v_t$ with $\gamma_t = \min\left\{\frac{\innp{\nabla f(x_t)}{x_t - v_t}}{M \norm{x_t - v_t}^2},1\right\}$ for some $M$. If $\innp{\nabla f(x_{t+1})}{x_t - v_t} \geq \frac{1}{2} \innp{\nabla f(x_{t})}{x_t - v_t}$, then it holds:
	\begin{equation*}
		f(x_t) - f(x_{t+1}) \geq \gamma_t \innp{\nabla f(x_t)}{x_t - v_t}/2 =
		\begin{cases*}
			\frac{\innp{\nabla f(x_t)}{x_t - v_t}^2}{2M \norm{x_t - v_t}^2} & if \quad $\gamma_t \in [0, 1)$ \\
			\innp{\nabla f(x_t)}{x_t - v_t} / 2 & if \quad $\gamma_t = 1$
		\end{cases*}.
	\end{equation*}
	\begin{proof}
		The proof follows directly via convexity and plugging in the definitions:
		\begin{align*}
			f(x_t) - f(x_{t+1}) & \geq \innp{\nabla f(x_{t+1})}{x_t - x_{t+1}} & \text{(convexity)} \\
			& \geq \gamma_t \innp{\nabla f(x_{t+1})}{x_t - v_t} & \text{(definition of $x_{t+1}$)} \\
			& \geq \gamma_t \innp{\nabla f(x_{t})}{x_t - v_t} / 2 & \text{(assumption of \eqref{eq:altAdaptive-simple})} \\
			& = \begin{cases*}
				\frac{\innp{\nabla f(x_t)}{x_t - v_t}^2}{2M \norm{x_t - v_t}^2} & if \quad $\gamma_t \in [0, 1)$ \\
				\innp{\nabla f(x_t)}{x_t - v_t} / 2 & if \quad $\gamma_t = 1$
			\end{cases*}. & \text{(definition of $\gamma_t$)}
		\end{align*}
	\end{proof}
\end{lemma}

Note that the proof above (again) explicitly relies on the convexity of $f$. It remains to show that \eqref{eq:altAdaptive-simple} holds for $2L$, whenever the function is $L$-smooth and $\gamma_t$ is the corresponding short step for $M = 2L$. The proof is very similar to before, however the last step is different.

\begin{lemma}[Smoothness implies \eqref{eq:altAdaptive-simple}] 
	\label{lem:smootToAdaptive-simple}	
	Let $f$ be $L$-smooth. Further, let $x_{t+1} = (1-\gamma_t) x_t + \gamma_t v_t$ with $\gamma_t = \min\left\{\frac{\innp{\nabla f(x_t)}{x_t - v_t}}{M \norm{x_t - v_t}^2},1\right\}$ and $M = 2L$. Then \eqref{eq:altAdaptive-simple} holds, i.e.,
	\begin{equation*}
		\innp{\nabla f(x_{t+1})}{x_t - v_t} \geq \frac{1}{2} \innp{\nabla f(x_{t})}{x_t - v_t}.	
	\end{equation*}
	\begin{proof}
		We use the alternative definition of smoothness using the gradients, i.e., we have 
		$$\innp{\nabla f(y) - \nabla f(x)}{y - x} \leq L \norm{y - x}^{2} 	\quad \text{for all } x, y \in P,$$
		by Remark~\ref{rem:smoothSC}. Now plug in $x \leftarrow x_t$ and $y \leftarrow x_{t+1}$, so that we obtain
		\begin{align*}
			\innp{\nabla f(x_{t+1}) - \nabla f(x_t)}{x_{t+1} - x_t} & \leq L \norm{x_{t+1} - x_t}^{2}
		\end{align*}
		and with plugging in the definition of $x_{t+1}$ we obtain
		\begin{align*}
			\innp{\nabla f(x_{t+1}) - \nabla f(x_t)}{\gamma_t(v_t - x_t)} & \leq L \gamma_t^2 \norm{v_t - x_t}^{2}.
		\end{align*}
		If $\gamma_{t} > 0$, dividing by $\gamma_t$ and then plugging in the definition of $\gamma_t$ yields
		$$
		\innp{\nabla f(x_{t+1}) - \nabla f(x_t)}{v_t - x_t} \leq \frac{1}{2} \innp{\nabla f(x_t)}{x_t - v_t},
		$$
		and rearranging gives the desired inequality 
		$$
		\innp{\nabla f(x_{t+1})}{x_t - v_t} \geq \frac{1}{2} \innp{\nabla f(x_t)}{x_t - v_t}.
		$$
		In case $\gamma_{t} = 0$ we have $x_t = x_{t+1}$ and the assertion holds trivially.
	\end{proof}
\end{lemma}

We will now show that \eqref{eq:altAdaptive-simple} is indeed weaker than \eqref{eq:altAdaptive} and that we cannot replace $M = L$ in Lemma~\ref{lem:smootToAdaptive}. To this end consider the following $1$-dimensional example: Pick $f(x) \doteq x^2$, so that $L=2$ holds. Consider $f: [-1, 1] \mapsto \R$ and $x_t = 1$. Then we have $\nabla f(x_t) = 2$, $v_t = -1$, and 
$$
\frac{\innp{\nabla f(x_t)}{x_t - v_t}}{L \norm{x_t - v_t}^2} = \frac{2(1 - (-1))}{2 (1 - (-1))^2} = \frac{1}{2},
$$
so that $\gamma_t = \min\left\{\frac{\innp{\nabla f(x_t)}{x_t - v_t}}{L \norm{x_t - v_t}^2},1\right\} = \frac{1}{2}$, $x_{t+1} = 0$, and $\nabla f(x_{t+1}) = 0$. This contradicts $0 = \innp{\nabla f(x_{t+1})}{x_t - v_t} \geq \frac{1}{2} \innp{\nabla f(x_t)}{x_t - v_t} > 0$.

% man hat x_{t+1} = 0, und nabla f(x_{t+1}) = 0.

\begin{algorithm}[h]
	\caption[]{(modified) Adaptive step-size strategy - simple variant}
	\label{alg:AdaptiveStepSize-simple}
	\begin{algorithmic}[1]
		\REQUIRE Objective function \(f\),
		smoothness estimate \(\widetilde{L}\),
		feasible points \(x\), \(v\)
		with \(\innp{\nabla f(x)}{x - v} \geq 0\),
		progress parameters \(\eta \leq 1 < \tau\)
		\ENSURE Updated estimate \(\widetilde{L}^{*}\),
		step-size \(\gamma\)
		\STATE\(M \leftarrow \eta \widetilde{L}\)
		\LOOP
		\STATE
		\(\gamma \leftarrow
		\min \{\innp{\nabla f(x)}{x - v}
		\mathbin{/} (M \norm{x - v}^{2}), 1 \}\) \COMMENT{compute short step for estimation $M$}
		%		\IF{\(f(x + \gamma (v - x)) - f(x)
			%			\leq \gamma \innp{\nabla f(x)}{v - x}
			%			+ \frac{\gamma^{2} M}{2} \norm{x - v}^{2}\)} 
		\IF{\(\innp{\nabla f(x + \gamma (v - x))}{x - v} \geq \frac{1}{2} \innp{\nabla f(x)}{x - v}\)} 
		\label{SufficientDecreaseSimple}
		\STATE \(\widetilde{L}^{*} \leftarrow M\)
		\RETURN \(\widetilde{L}^{*}\), \(\gamma\)
		\ENDIF
		\STATE\(M \leftarrow \tau M\)
		\ENDLOOP
	\end{algorithmic}
\end{algorithm}

\end{document}

%% file: prologue.tex
\usepackage[utf8]{inputenc} % allow utf-8 input
\usepackage[T1]{fontenc}    % use 8-bit T1 fonts

 \usepackage[dvipsnames]{xcolor}
\usepackage[
	bookmarks=true,
	bookmarksopen=true,
	bookmarksopenlevel=3,
	bookmarksnumbered=true,
	plainpages=false,
	colorlinks=true,
	urlcolor=DarkBlue,
	citecolor=DarkBlue,
	linkcolor=DarkBlue
]{hyperref}

\usepackage{booktabs}       % professional-quality tables
\usepackage{amsfonts}       % blackboard math symbols
\usepackage[protrusion=true,expansion,tracking=true]{microtype}     % microtype
\usepackage{nicefrac} 
\usepackage{amsmath}
\usepackage{amsthm}
\usepackage{amssymb}
\usepackage{thmtools}
\usepackage{paralist}
\usepackage{comment}
\usepackage{dsfont}

\usepackage{bbm}
\usepackage{float}
\usepackage{balance}
\usepackage{stmaryrd}

\usepackage{multicol}
\usepackage[export]{adjustbox}
\usepackage{multirow}
\usepackage{tabularx}
\usepackage{mathtools}
\usepackage{mdframed}

\usepackage[shortlabels]{enumitem}
\setitemize{noitemsep,topsep=0pt,parsep=0pt,partopsep=0pt}

\usepackage[ruled, linesnumbered]{algorithm2e} %
\usepackage{algorithmic}
\algsetup{
	linenosize=\scriptsize,
	linenodelimiter={~}
}

\DontPrintSemicolon

\usepackage{jmlr2e}

\usepackage[margin=1in]{geometry}

\usepackage[capitalize,nameinlink]{cleveref}[0.21]

\usepackage{natbib}
\setcitestyle{authoryear,round,citesep={;},aysep={,},yysep={;}}
\renewcommand{\cite}[1]{\citep{#1}}

\usepackage{libertine}
\usepackage{libertinust1math}
\usepackage{inconsolata}

\RequirePackage{pifont}
\RequirePackage{dsfont} % or bbold redefining \mathbb

\ifedit
\newcommand{\pagehint}[1]{\quad\textcolor{blue}{\footnotesize\textup{[#1 page(s)]}}}
\newcommand{\hint}[1]{\quad\textcolor{blue}{\footnotesize\textup{[#1]}}}
\usepackage{showlabels}
\else
\newcommand{\pagehint}[1]{}
\newcommand{\hint}[1]{}
\fi

\usepackage{graphicx,wrapfig,lipsum}
\usepackage[
singlelinecheck=false, % <-- important
labelfont=bf
]{caption}
\usepackage{makecell}
\usepackage[dvipsnames]{xcolor}

\definecolor{DarkMagenta}{rgb}{0.55, 0, 0.55}
\definecolor{DarkBlue}{rgb}{0, 0, 0.55}

\def\RR{{\mathbb{R}}}

\renewcommand{\H}{\mathbb{H}}
\newcommand{\N}{\mathbb{N}}

\newcommand{\R}{\mathbb{R}}

\DeclareMathOperator{\argmax}{argmax}
\DeclareMathOperator{\argmin}{argmin}
\DeclareMathOperator{\conv}{conv}

\DeclarePairedDelimiterXPP{\innp}[2]{\mathinner\bgroup}{\langle}{\rangle}%
{\egroup}{#1, #2}

\newcommand{\norm}[1]{\| #1 \|}

\newcommand{\card}[1]{| #1 |}

\theoremstyle{plain} \numberwithin{equation}{section}
\newtheorem{theorem}{Theorem}[section]
\numberwithin{theorem}{section}
\newtheorem{lemma}[theorem]{Lemma}

\newtheorem{proposition}[theorem]{Proposition}

\theoremstyle{definition}
\newtheorem{definition}[theorem]{Definition}

\newtheorem{remark}[theorem]{Remark}

\newcommand{\newalgorithm}[2]{%
	\newenvironment{#1}[1][]{%
		\begin{algorithm}[{##1}]\captionsetup{name={#2},list=no}%
		}{\end{algorithm}}}

\newalgorithm{oracle}{Oracle}

\theoremstyle{plain}